\documentclass{gtart_h}

\def\ifplaintex{\expandafter\ifx\csname documentclass\endcsname\relax}

\def\gtp{{\mathsurround=0pt\it $\cal G\mskip-2mu$eometry \&\ 
$\cal T\!\!$opology $\cal P\!$ublications}}  

\def\recd{{\small Received:\qua\receiveddate\ifx\reviseddate\relax
\else\qquad Revised:\qua\reviseddate\fi\par}} 

\def\lognumber#1{\def\thelognumber{#1}}
\def\volumenumber#1{\def\thevolumenumber{#1}}
\def\volumeyear#1{\def\thevolumeyear{#1}}
\def\papernumber#1{\def\thepapernumber{#1}}
\def\pagenumbers#1#2{\def\startpage{#1}\def\finishpage{#2}}
\def\published#1{\def\publishdate{#1}}

\def\received#1{\def\receiveddate{#1}}

\def\accepted#1{\def\accepteddate{#1}}

\def\asciiaddress#1{\def\theasciiaddress{#1}}
\def\asciiemail#1{\def\theasciiemail{#1}}

\let\\\par\let\thelognumber\relax\let\thevolumenumber\relax
\let\thepapernumber\relax\let\thevolumeyear\relax\let\startpage\relax
\let\finishpage\relax\let\publishdate\relax\let\receiveddate\relax
\let\reviseddate\relax\let\accepteddate\relax\let\theasciititle\relax
\let\theasciiauthors\relax\let\theasciiaddress\relax
\let\theasciiabstract\relax

\let\theasciiemail\relax

\ifplaintex
\font\logobig=cmssbx10 scaled 3836
\font\logomed=cmssbx10 scaled 2557
\else
\font\logobig=cmssbx10 scaled 4200
\font\logomed=cmssbx10 scaled 2800
\fi

\long\def\makeagttitle{   
\count0=\startpage
\agt\hfill      
\hbox to 45truept{\vbox to 0pt{\vglue -13truept{\logomed A\kern -.37em{\logobig 
T}\kern -.38em G}\vss}\hss}
\break
{\small Volume \thevolumenumber\ (\thevolumeyear)
\startpage--\finishpage\nl
Published: \publishdate}

\vglue .25truein

{\parskip=0pt\leftskip 0pt plus
1fil\def\\{\par\smallskip}{\Large\bf\thetitle}\par\medskip} \vglue
0.05truein

{\parskip=0pt\leftskip 0pt plus 1fil\def\\{\par}{\sc\theauthors}
\par\medskip}%
 
\vglue 0.03truein 

{\small\leftskip 25truept\rightskip 25truept{\bf Abstract}\stdspace\theabstract

{\bf AMS Classification}\stdspace\theprimaryclass
\ifx\thesecondaryclass\relax\else; \thesecondaryclass\fi\par
{\bf Keywords}\stdspace \thekeywords\par}\vglue 7truept

}   

\ifplaintex
\font\phead=cmsl9 scaled 950
\font\pnum=cmbx10 scaled 913
\font\pfoot=cmsl9 scaled 950
\headline{\vbox to 0pt{\vskip -4.5mm\line{\small\phead\ifnum
\count0=\startpage ISSN 1472-2739 (on-line) 1472-2747 (printed)
\hfill {\pnum\folio}\else\ifodd\count0\def\\{ }%
\ifx\theshorttitle\relax\thetitle\else\theshorttitle\fi\hfill{\pnum\folio}
\else\def\\{ and }{\pnum\folio}\hfill\ifx\theshortauthors\relax\theauthors
\else\theshortauthors\fi\fi\fi}\vss}}
\footline{\vbox to 0pt{\vglue 0mm\line{\small\pfoot\ifnum\count0=\startpage
\copyright\ \gtp\hfill\else
\agt, Volume \thevolumenumber\ (\thevolumeyear)\hfill\fi}\vss}}
\else
\font=cmsl9 scaled 1050
\font\lnum=cmbx10 
\font=cmsl9 scaled 1050
\makeatletter
\def\@oddhead{{\small\lhead\ifnum\count0=\startpage ISSN 1472-2739 
(on-line) 1472-2747 (printed)\hfill {\lnum\number\count0}\else\ifodd\count0
\def\\{ }\ifx\theshorttitle\relax \thetitle \else\theshorttitle\fi\hfill
{\lnum\number\count0}\else\def\\{ and }{\lnum\number\count0}
\hfill\ifx\theshortauthors\relax 
\theauthors\else\theshortauthors\fi\fi\fi}}\def\@evenhead{\@oddhead}
\def\@oddfoot{\small\lfoot\ifnum\count0=\startpage\copyright\ \gtp\hfill\else
\agt, Volume \thevolumenumber\ (\thevolumeyear)\hfill\fi}
\def\@evenfoot{\@oddfoot}
\makeatother
\fi
\let\maketitlepage\makeagttitle
\let\makeshorttitle\maketitlepage
\let\maketitle\maketitlepage

\newwrite\gtoutfile
\long\gdef\makeheadfile{  
{\def\\{, }\def\s{ }
\immediate\openout\gtoutfile head.xxx
\immediate\write\gtoutfile{Proxy-for: \ifx\theasciiauthors\relax
\theauthors\else\theasciiauthors\fi\s<\ifx\theasciiemail\relax\theemail\else\theasciiemail\fi>}
\immediate\write\gtoutfile{\noexpand\\}
\immediate\write\gtoutfile{Authors: \ifx\theasciiauthors\relax
\theauthors\else\theasciiauthors\fi}
{\def\\{ }\immediate\write\gtoutfile{Title: \ifx\theasciititle\relax
\thetitle\else\theasciititle\fi}}
\immediate\write\gtoutfile{Subj-class: GT or SG, GR etc}
\immediate\write\gtoutfile{MSC-class: \theprimaryclass\ifx\thesecondaryclass\relax\else, \thesecondaryclass\fi}
\immediate\write\gtoutfile{Journal-ref: Algebr. Geom. Topol. \thevolumenumber\s
(\thevolumeyear) \startpage-\finishpage}
\immediate\write\gtoutfile{Comments: Published by Algebraic and
Geometric Topology at}
\immediate\write\gtoutfile{\s\s\s  http://www.maths.warwick.ac.uk/agt/AGTVol\thevolumenumber/agt-\thevolumenumber-\thepapernumber.abs.html}
\immediate\write\gtoutfile{\noexpand\\}
\immediate\write\gtoutfile{}
\ifx\theasciiabstract\relax
\immediate\write\gtoutfile{\theabstract}\else
\immediate\write\gtoutfile{\theasciiabstract}\fi
\immediate\write\gtoutfile{}
\immediate\write\gtoutfile{\noexpand\\}
\immediate\write\gtoutfile{}
\immediate\closeout\gtoutfile}}  

\def\maketitlepage{\makeagttitle\makeheadfile}
\let\makeshorttitle\maketitlepage
\let\maketitle\maketitlepage

\lognumber{21}
\volumenumber{4}
\volumeyear{2004}
\papernumber{21}
\published{24 June 2004}
\pagenumbers{407}{437}
\received{8 June 2004}
\accepted{10 June 2004}

\usepackage{amsmath,amssymb,epsf}

\input xy
\xyoption{all}

\newtheorem{thm}{Theorem}[section]
\newtheorem{lem}[thm]{Lemma}

\newtheorem{prop}[thm]{Proposition}

\theoremstyle{definition}
\newtheorem{defn}[thm]{Definition}

\def\bZ{\mathbb{Z}}
\def\bQ{\mathbb{Q}}

\def\bC{\mathbb{C}}

\def\cL{\mathcal{L}}

\def\cC{\mathcal{C}}
\def\cS{\mathcal{S}}
\def\cD{\mathcal{D}}

\def\cE{\mathcal{E}}

\def\cT{\mathcal{T}}

\newcommand{\sn}{S_n}
\newcommand{\sm}{S_m}
\newcommand{\cn}{C_n}
\newcommand{\cm}{C_m}

\newcommand{\op}[1]{\overline{#1}}
\newcommand{\sga}{\Sigma}

\newcommand{\ccc}{\colon}
\newcommand{\ra}{\rightarrow}

\newcommand{\gSXone}{g_1\ccc\sga_1\ra X}
\newcommand{\gSXtwo}{g_2\ccc\sga_2\ra X}

\newcommand{\thsx}{\cT_X}

\newcommand{\vect}{\text{Vect}}
\newcommand{\ehqft}{E\ccc\thsx\ra\vect_{1}}

\newcommand{\lan}{\langle}
\newcommand{\ran}{\rangle}

\newcounter{test}

\DeclareMathOperator{\Id}{Id}
\DeclareMathOperator{\MSO}{MSO}
\DeclareMathOperator{\Hom}{Hom}

\DeclareMathOperator{\Ext}{Ext}
\DeclareMathOperator{\Ker}{Ker}
\DeclareMathOperator{\Tift}{TIFT}
\newcommand{\vpic}[2]{\epsfxsize #2\epsffile{#1.eps}}
\newcommand{\Ktor}{\{\text{$K^\times$--torsors}\}}
\newcommand{\Kmod}{\{\text{rank one $K$--modules}\}}

\begin{document}
\title{Gerbes and homotopy quantum field theories}
\author{Ulrich Bunke\\Paul Turner\\Simon Willerton}

\asciiemail{bunke@uni-math.gwdg.de, paul@ma.hw.ac.uk, s.willerton@shef.ac.uk}  
\gtemail{\mailto{bunke@uni-math.gwdg.de}, \mailto{paul@ma.hw.ac.uk},
\mailto{s.willerton@shef.ac.uk}}

\asciiaddress{Mathematisches Institut, Universitaet Goettingen, 
37073 Goettingen, Germany\\Department of Mathematics, 
Heriot-Watt University, Edinburgh EH14 4AS, Scotland\\Department 
of Pure Mathematics, University of Sheffield, Sheffield S3 7RH, England}
\address{Mathematisches Institut, Universit\"at G\"ottingen, 
37073 G\"ottingen, Germany\\\smallskip\\Department of Mathematics, 
Heriot-Watt University, Edinburgh EH14 4AS, Scotland\\\smallskip\\Department 
of Pure Mathematics, University of Sheffield, Sheffield S3 7RH, England}
\shortauthors{Ulrich Bunke, Paul Turner and Simon Willerton}


\begin{abstract}
For smooth finite dimensional manifolds, we characterise gerbes with
connection as functors on a certain surface cobordism category.  This
allows us to relate gerbes with connection to Turaev's 1+1-dimensional
homotopy quantum field theories, and we show that flat gerbes are
related to a specific class of rank one homotopy quantum field
theories.
\end{abstract}
\maketitle\let\\\par

\section*{Introduction}
\addcontentsline{toc}{section}{Introduction}

The original motivation for this paper was to reconcile the two
``higher'' versions of a line bundle with connection mentioned in
the title. In the process we came up with a characterization of
gerbes-with-connection over a fixed space as functors from a
certain cobordism category.  Before getting onto that we will give
a quick description of the two objects in the title, but first it
is pertinent to give a reminder of what a line bundle with
connection is.


\subsection*{Line bundles}
\addcontentsline{toc}{subsection}{Line bundles}
A line bundle with connection can be viewed in many ways, especially
if one wants to generalise to ``higher'' versions.  Here we will
mention the idea of it being determined by its holonomy, of it being a
functor on the path category of the base space and of it being a
functor on the $0+1$ dimensional cobordism category on the base.

The holonomy of a line bundle with connection is a $\bC^\times$-valued
function on the free loop space of the base manifold $X$.  Barrett
\cite{Barrett:HolonomyAndPathStructures} and others showed that
functions on free loop space which occur as the holonomy of a line
bundle with connection are characterised by being invariant under
thin-homotopy (see below), being invariant under diffeomorphism of the
circle, and satisfying a smoothness condition.  Further, such a
function on the free loop space uniquely determines a line bundle with
connection up to equivalence.

A different characterisation of holonomy is given by thinking of it as
a function $\widetilde S$ on smooth maps of closed one-dimensional
manifolds into the space which is multiplicative under disjoint union,
together with a closed two-form $\widetilde c$, the curvature of the
bundle, such that if $\widetilde v\colon\widetilde V\to X$ is a map of
a surface with boundary into the base space then $\widetilde
S(\partial \widetilde v)=\exp (\int_{\widetilde V} \widetilde
v^*\widetilde c)$.  Such pairs $\widetilde S$ and $\widetilde c$
satisfying this condition form the Cheeger-Simons group (see
\cite{cheegersimons83}) of differential characters $\widehat H^2(X)$.
This group parametrises line bundles with connection up to
equivalence.

One could take a groupoid point of view of a line bundle with
connection in the following fashion.  The path category $PX$ of a
space $X$ is the category whose objects are the points in the space
and whose morphisms are, roughly speaking, smooth paths between them,
while the category $\text{Vect}_{1}$ is the category of
one-dimensional complex vector spaces with invertible linear maps as
morphisms.  Any line bundle with connection gives a functor $PX\to
\text{Vect}_{1}$, which to a point in $X$ associates the fibre over
that point, and to a path between two points associates the parallel
transport along that path.  This functor will satisfy some smoothness
condition and will in fact descend to a functor on the thin-homotopy
path groupoid.  Actually, here the category $\text{Vect}_{1}$ is
rather large and could be replaced by something like the small
category of lines in infinite projective space.

A variation on this is obtained by considering the 0+1-dimensional
cobordism category of the space $X$, this has finite ordered
collections of points in $X$ as its objects, and cobordisms between
them as morphisms.  A monoidal functor from this to the category of
complex lines, with the usual tensor product, should be an appropriate
notion of bundle with connection.

Things become a lot simpler when\/ {\em flat\/} bundles are
considered.  In this case geometric notions descend to topological
ones.  A flat line bundle is a line bundle with connection whose
curvature vanishes identically.  This means that the holonomy can be
considered as an element of the first cohomology group
$H^1(X,\bC^\times)$.  The categorical descriptions become a lot
simpler because the morphism sets can be quotiented by homotopy
relations, thus becoming discrete sets.


\subsection*{Gerbes}
\addcontentsline{toc}{subsection}{Gerbes}
A gerbe is essentially a realization of a degree three cohomology
class.  The idea of a gerbe was brought to many people's attention
by the book of Brylinski \cite{Brylinski:LoopSpacesBook}.  There
are several different but equivalent ways to realize gerbes, these
include ways involving sheaves of categories
\cite{Brylinski:LoopSpacesBook}, classifying space bundles
\cite{Gajer:GeometryOfDeligneCohomology}, bundle-gerbes
\cite{Murray:BundleGerbes}, and bundle realizations of Cech
cocycles
\cite{Chaterjee:Thesis,Hitchin:SpecialLagrangianSubmanifolds}. The
degree three cohomology class corresponding to a gerbe is called
its Dixmier-Douady class and it is the analogue of the first Chern
class for a line bundle.

It is possible to introduce an appropriate notion of connection (by
which we mean a curving and a connective structure in the language of
\cite{Brylinski:LoopSpacesBook}) on a gerbe over a smooth manifold.
Associated to a gerbe with connection, in any of the above mentioned
descriptions, are a curvature three-form and two notions of holonomy.
The first notion of holonomy is a $\bC^\times$-valued function $S$ on
the space of maps of closed surfaces into the manifold; and the second
is loop holonomy which is a line bundle with connection on the free
loop space of the manifold.  The curvature and surface holonomy are
related by the fact that if $v\colon V\to X$ is a map of a three
manifold with boundary into the base manifold then $S(\partial v)=\exp
(2\pi i \int_V v^*c)$.  Again, pairs $S$ and $c$ which satisfy this
condition form a group $\widehat H^3(X)$ which parametrises gerbes
with connection up to equivalence.  (This group is isomorphic to the
so called Cheeger-Simons group of differential characters, see
Appendix \ref{ap:csg} for more details.)  In this paper we will
work with this holonomy description of gerbes.

A\/ {\em flat gerbe\/} is a gerbe with connection whose curvature
vanishes identically.  This implies that the Dixmier-Douady class of the
underlying gerbe is a torsion class.  It  also follows  that the
holonomy around a surface only depends on the homology class of the
surface, and so the holonomy can be considered as an element of
$\operatorname{Hom}(H_2(X,\bZ), \bC^\times) \cong H^2(X,\bC^\times)$.
This establishes a bijection between flat gerbes on $X$ and
$H^2(X,\bC^\times)$.


\subsection*{Homotopy quantum field theories}
\addcontentsline{toc}{subsection}{Homotopy quantum field theories}
The second generalisation is the notion of a 1+1-dimensional homotopy
quantum field theory, which strictly speaking generalises the idea of
a\/ {\em flat\/} bundle.  This notion was introduced by Turaev in
\cite{Turaev:HomotopyFieldTheoryInDimension2} and independently by
Brightwell and Turner
\cite{BrightwellTurner:RepresentationOfHomotopySurfaceCategory}, but
the idea goes back to the work of Segal
\cite{Segal:EllipticCohomology}.  Turaev considered the case of
Eilenberg-MacLane spaces and, orthogonally, Brightwell and Turner
considered simply-connected spaces.  A homotopy
quantum field theory is like a topological quantum field theory taking
place in a ``background space'', and it can be
given a functorial description as follows.\footnote%
{It should be noted that we alter Turaev's definition by removing Axiom
2.27 which is not appropriate for non-Eilenberg-MacLane spaces, this
does not alter any of the theorems in his paper provided they are all
stated for Eilenberg-MacLane spaces. This is the position adopted
by Rodrigues in \cite{ro:hqft}, where a connection with gerbes and
thinness is also suggested.}
For a space, the homotopy surface category of a space generalizes the
$0+1$-dimensional cobordism category by having as objects collection
of loops in the space and having as morphisms cobordisms between these
considered up to boundary preserving homotopy.  A 1+1-dimensional
homotopy quantum field theory on a space is a symmetric monoidal functor from
the homotopy surface category of a space to the category of finite
dimensional vector spaces.  As we are comparing these to gerbes we
will only consider the rank one theories, ie.\ those functors taking
values in the subcategory of one dimensional vector spaces.


\subsection*{Outline of this paper}
\addcontentsline{toc}{subsection}{Outline of this paper}
One motivation for this paper was to figure out how gerbes are
related to homotopy quantum field theories, another was to
understand what conditions are necessary and sufficient for a line
bundle on loop space to come from a gerbe.  These questions are
addressed by considering an object that we have dubbed a
thin-invariant field theory.  The main novelty is that it uses the
idea of thin-cobordism: two manifolds in a space are thin
cobordant if there is a cobordism between them which has ``zero
volume'' in the ambient space.  In $1$-dimensional manifolds this
is the same as thin-homotopy as defined by Barrett
\cite{Barrett:HolonomyAndPathStructures} and further developed in
\cite{CaetanoPicken:AxiomaticHolonomy} and
\cite{MackaayPicken:HolonomyAndParallelTransport}, but this is not
the case for surfaces --- the 1+1-dimensional thin-cobordism
category is a groupoid whereas the 1+1-dimensional thin-homotopy
category is not.  A thin-invariant field theory is essentially a
smooth symmetric monoidal functor from the thin-cobordism category
to the category of one-dimensional vector spaces.  The idea is
that this gives an alternative description of a gerbe.  The view
that a gerbe should be a functor on a cobordism category has been
advocated by Segal in \cite{Segal:TopologicalStructures}. The
collection of thin-invariant field theories on $X$ form a group in
a natural way.  Our main theorem is the following.
\setcounter{section}{6} \setcounter{thm}{3}\addtocounter{thm}{-1}
\begin{thm}
 On a smooth finite dimensional manifold, there is
  an isomorphism from the group of thin-invariant field theories (up to
  equivalence) to the group of gerbes with connection (up to equivalence).
\end{thm}

To make the connection with homotopy quantum field theories we
show that a certain natural subset, the rank one, normalised ones,
correspond to flat thin-invariant field theories, this gives the following.
\setcounter{section}{6}
\setcounter{thm}{4}\addtocounter{thm}{-1}
\begin{thm}
  On a finite dimensional manifold, the group of normalised rank one homotopy
  quantum field theories (up to equivalence) is isomorphic to the
  group of flat gerbes (up to equivalence).
\end{thm}
In this context it makes sense to consider homotopy quantum field
theories defined over an arbitrary commutative ring with unity.  We
classify these in the following manner, which generalizes a theorem of
Turaev.
\setcounter{section}{7}
\setcounter{thm}{1}\addtocounter{thm}{-1}
\begin{thm}
  Let $K$ be a commutative ring, and $X$ be a path
  connected topological space.  Then Turaev's construction gives an
  isomorphism between the group $H^2(X,K^\times)$ and the group of
  normalised, rank one homotopy quantum field theories defined over
  $K$.
\end{thm}

We then show that a thin-invariant field theory  is an extension of
the usual line bundle of
a gerbe over the free loop space.
\setcounter{section}{8}
\setcounter{thm}{1}\addtocounter{thm}{-1}
\begin{thm}
A thin-invariant field theory  can be restricted to the path category
of the free loop
space giving  a line bundle with connection on the free loop
space. This is isomorphic to the  transgression of the associated gerbe.
\end{thm}

We include two appendices. In the first we compare our definition of
the Cheeger-Simons group with the more familiar one, and in the second
we gather together, for ease of reference, a number of categorical
definitions used throughout the paper.

It is worth noting here that homotopy quantum field theories are
bordism-like in their nature, whereas gerbes are homological
creatures.  It seems to us that the techniques used and results
obtained in this paper rely on the coincidence of bordism and homology
in low degree, and will not necessarily generalise to higher degrees.

\setcounter{section}{0}


\section{Basic definitions}
\subsection{Bordism}
Here we give, for those unfamiliar with the notion, a very brief
introduction to (co)bordism groups and then we present the
low-dimensional co-incidence result which is central to the paper.

The $n$th oriented bordism%
\footnote{There is a standard problem with terminology here.
Initially bordism groups were called cobordism groups, because two
things are cobordant if they cobound something else. Unfortunately
in this context the prefix ``co'' usually refers to the
contravariant theory, so cobordism was taken to mean the
contravariant version (analogous to cohomology) and the word
bordism was used for the covariant theory (analogous to homology).
In this paper we will always be interested in the covariant
theory.}\ %
 group $\MSO_n(X)$ of a space $X$ for a
non-negative integer $n$ is similar but subtly different to the
$n$th ordinary homology group $H_n(X)$.  Whereas homology groups
are defined using chains of simplices, bordism groups are defined
using maps of manifolds. The main ingredient in the definition is
the set of pairs $(V,v)$ where $V$ is an oriented smooth
$n$--manifold and $v\colon V\to X$ is a map.  Two such pairs
$(V,v)$ and $(V',v')$ are said to be \emph{cobordant} if there is
an $(n+1)$--manifold $W$ with $\partial W\cong \op{V}\sqcup V'$
 and a map $w\colon W\ra X$ such that
$\partial w = \op{v} \sqcup v^\prime $.  The group $MSO_n(X)$ is defined to be the set
of equivalence classes under this cobordism relation, with the group structure being
induced from the disjoint union of manifolds.

These groups share many properties with ordinary homology groups
$H_n(X)$, forming an example of what is called an extraordinary
homology theory.  In fact the only difference between homology and
bordism lies in the torsion part, as rationally they are the same:
$H_n(X)\otimes \bQ\cong \MSO_n(X)\otimes \bQ$.  The general theory
is a well developed topic in the algebraic topology literature,
one source for a comprehensive treatment would be \cite{switzer}.

 The following lemma on the low-dimensional
co-incidence of bordism and homology is key to the ideas of this
paper.

\begin{lem}
  \label{Lemma:SecondHomologyEqualsBordism}
  The first homology and bordism groups of a space are isomorphic, as
  are the second groups: if $X$ is a space then
  $\MSO_1(X)\cong H_1(X;\bZ)$ and $\MSO_2(X)\cong H_2(X;\bZ)$.
\end{lem}
\begin{proof}
  Apply the Atiyah-Hirzebruch spectral sequence for
  bordism groups (see for example \cite{switzer})
  and use the fact that in low dimensions the
  coefficients for
  bordism are given by $\MSO_0(\text{pt})\cong \bZ$ and
  $\MSO_1(\text{pt})\cong\MSO_2(\text{pt})\cong\{1\}$.
\end{proof}

\subsection{$X$-surfaces and thin cobordisms} In this section we
introduce the key notions of $X$-surfaces, thin-cobordism and thin
homotopy.

If $X$ is a smooth manifold then an\/ {\em $X$-surface\/} is
essentially a smooth map of a surface into $X$, but with certain
technical collaring requirements to ensure that $X$-surfaces can
be glued together.  It is perhaps possible to avoid these
technical conditions by working with piece-wise smooth maps, but
we have not done that.

Boundaries of surfaces will need to be parametrised, so for
concreteness, let $S^1$ be the set of unit complex numbers and fix an
orientation for this.  Let $S_n$ be the union of $n$ ordered copies of
$S^1$.  Fix an orientation on the unit interval $[0,1]$ and define
$C_n:=S_n\times [0,1]$, so that $C_n$ is $n$ ordered parametrised
cylinders.  For any oriented manifold $Y$ let $\overline Y$ denote the
same manifold with the opposite orientation.

A surface will mean a smooth oriented two-manifold $\Sigma$ together
with a\/ {\em collar\/}, which will mean a certain type of
parametrisation of a neighbourhood $N_\Sigma$ of the boundary:
\[
   \iota_{\sga}\ccc\op \cm \sqcup \cn
            \stackrel\cong\longrightarrow N_\Sigma.
\]
The $m$ boundary components corresponding to $\op \cm$ will be
called\/ {\em inputs\/} and the $n$ corresponding to $\cn$ will be
called\/ {\em outputs}. Note that inputs and outputs inherit an order
from $\iota_{\sga}$.

Define an $X$-surface to be a surface $\Sigma$ as above, and a smooth
map $g\colon \Sigma \to X$ such that $g|_{N_\Sigma} \circ
\iota_\Sigma$ factors through the projection $\op \cm \sqcup \cn \to
\op \sm \sqcup \sn$, ie.\ the map $g$ is constant in transverse
directions near the boundary.  The inputs and outputs of $g$ are the
restrictions of $g$ to the inputs and outputs of the underlying
surface.  If the inputs of $\gSXone$ agree with the outputs of
$\gSXtwo$ then we can glue $\Sigma_1$ and $\Sigma_2$ together using the
given collars to form another surface $\Sigma_1 \cup \Sigma_2$ and
using the induced maps form the $X$-surface $g_1 \circ g_2\ccc \Sigma_1
\cup \Sigma_2 \ra X$.  If the inputs of $g$ are the same as the
outputs, then use the notation $\lan g\ran$ to denote the closed
$X$-surface obtained by gluing the inputs to the outputs.  An
$X$-three-manifold is defined similarly.

Informally two $X$-surfaces are thin cobordant if there exists a
cobounding manifold which has no volume in $X$. More formally,
  two $X$-surfaces $g\colon \Sigma\to X$ and $g'\colon
  \Sigma'\to X$ are {\em thin cobordant} if there exists a collared
  three-manifold $W$  such that $\partial W\cong \Sigma\cup
  \overline{\Sigma'}$ and a smooth map $w\colon W\to X$ satisfying
  $w|_{\partial W}=g\cup \overline{g'}$ and $dw$ everywhere having
  rank at most two.

Thin homotopy is a particular kind of thin cobordism.  Let
$g\colon\Sigma \to X$ and
$g'\colon\Sigma'\to X$ be $X$-surfaces with the same inputs and the
same outputs.  The
maps $g$ and $g^\prime$ are {\em thin homotopic} if there exists a
thin cobordism homotopic to $\Sigma\times[0,1]$.

One fundamental difference between thin-homotopy and
thin-cobordism is that cobordisms are invertible modulo
thin-cobordism, but not modulo thin-homo\-topy.  The next
proposition shows that if $g$ is a cobordism then its reversal
$\op g$ is an inverse modulo thin-cobordism.
\begin{prop}\label{prop:holonomycancel}
  If $g\colon \Sigma\to X$ is an $X$-surface which is not necessarily
  closed, then the closed $X$-surface $\lan g\circ \op g\ran$ is
  thin-cobordant to the empty $X$-surface.
\end{prop}
\begin{proof}
  Consider the manifold with corners $\Sigma\times I$.  Smooth
  this by just removing an arbitrarily small neighbourhood of the
  corners and call the resulting smooth manifold $W$.  The collaring
  implies that the boundary of
  $W$ can be identified with $\Sigma\cup \overline \Sigma$.  Define
  the map $w\colon W\to X$ to be the projection to $\Sigma$ composed
  with $g$.  The differential $dw$ automatically has rank at most two
  and thus $w$ provides the requisite thin cobordism.
\end{proof}


\section{Gerbe holonomy}
\label{Section:GerbeHolonomy}

In this section we collect together the facts we need about gerbe
holonomy.

For a gerbe with connection on a manifold $X$ there is the associated
gerbe holonomy which associates a complex number to each closed $X$-surface.
The gerbe holonomy is invariant under diffeomorphism of $X$-surfaces
and it is multiplicative under disjoint union.

The holonomy is related to the curvature of the gerbe connection in
the following fashion.  Suppose that $S$ is the gerbe holonomy and
$c$, a closed three-form, is the gerbe curvature.  If $v\colon V\to X$
is an $X$-three-manifold then the following holonomy-curvature
relation holds:
\[
  S(\partial v)= \exp \left(2\pi i\int_V v^*c\right).
\]

We can take all of the diffeomorphism invariant, multiplicative
functions on the set of closed $X$-surfaces for which there exists a
three-form so that the holonomy-curvature relation is satisfied.
These form a group $\widehat H^3(X)$.  This is not exactly the third
Cheeger-Simons group, which is defined using smooth two-cycles rather
than closed $X$-surfaces.  However these two groups are isomorphic in
this degree, this is proved in Appendix~\ref{ap:csg} and is due to the
fact the bordism and homology agree at low degree.  We will therefore
refer to $\widehat H^3(X)$ as the Cheeger-Simons group.  Thus each gerbe
with connection gives rise to an element in this Cheeger-Simons group
by means of its surface holonomy.  It turns out that this sets up a
bijection between gerbes with connection and this Cheeger-Simons group
(see eg.\ \cite{Gajer:GeometryOfDeligneCohomology}).  Thus specifying
a gerbe with connection is the same as specifying its surface
holonomy.  We will think of $\widehat H^3(X)$ as the group of gerbes
with connection.

There are two useful exact sequences involving gerbes wich we will now
mention (see \cite[Section~1.5]{Brylinski:LoopSpacesBook}).
Let $\Omega^*(X)$ denote the smooth complex differential forms on $X$.
By $\Omega^2(X)_{d=0, \bZ}$ we denote the subspace of closed forms
which have periods in $\bZ$.  There are the following exact sequences.
\begin{align*}
0\rightarrow H^2(X,\bC^\times)&
  \stackrel{\phi}{\rightarrow} \widehat H^3(X)
  \stackrel{c}{\rightarrow} \Omega^3(X);\\
0\rightarrow \Omega^2(X)/\Omega^2(X)_{d=0, \bZ}&
  \stackrel{h}{\rightarrow}\widehat H^3(X)
  \stackrel{D} \rightarrow H^3(X,\bZ)\rightarrow 0.
\end{align*}
Here $c$ associates to each gerbe its curvature, $D$ maps a gerbe to
its Dixmier-Douady class, and $h$ maps a class $[\omega]$, $\omega\in
\Omega^2(X)$, to the gerbe with curvature $d\omega$ and holonomy
$h([\omega])\lan g\ran=\exp(2\pi i \int_{\Sigma} g^* \omega)$ for all
$g:\Sigma\rightarrow X$ with $\Sigma$ an oriented closed surface.
The map $\phi$ can be interpreted as the inclusion of flat gerbes.

Mackaay and Picken \cite{MackaayPicken:HolonomyAndParallelTransport}
observed that gerbe holonomy is invariant under thin-homotopy: we make
the stronger, key observation that it is invariant under
thin-cobordism.
\begin{prop}\label{prop:holonomythininvariance}
  Suppose that $S$ is the holonomy of a gerbe with connection on a
  manifold $X$.  If $g\colon \Sigma\to X$ and $g'\colon \Sigma'\to X$
  are closed $X$-surfaces which are thin cobordant then $S(g)=S(g')$.
\end{prop}
\begin{proof}
  Let the three-form $c$ be the curvature of the gerbe and suppose
  that $w\colon W\to X$ is a thin cobordism between $g$ and $g'$.
  The holonomy-curvature relation implies that
  $S(\partial w)=\exp(2\pi i\int_W w^*c)$.  However, the right-hand side is
  equal to one as $dw$ has rank at most two, and the left-hand side is
  equal to $S(g\cup \overline{g'})= S(g)S(g')^{-1}$, from which the
  result follows.
\end{proof}


\section{Thin-invariant field theories}

In order to define thin-invariant field theories we adopt a similar
philosophy to
\cite{BrightwellTurner:RepresentationOfHomotopySurfaceCategory} (see
also \cite{ro:hqft}) and define a category of cobordisms in a background
$X$ and then define a thin-invariant field theory to be a
complex representation of this category.

\begin{defn}
  The {\em thin-homotopy surface category}, $\thsx$, of a smooth
  manifold $X$, is the category whose objects are smooth maps from $n$
  copies of a standard circle, $S_n$, to $X$ for some integer $n$; and
  for which a morphism is an $X$-surface $g\colon \Sigma \to X$
  considered as a morphism from the map defined by its inputs to the
  map defined by its outputs, with $g'\colon \Sigma' \to X$ being
  identified with $g$ if there exists a diffeomorphism $T\ccc \Sigma
  \ra \Sigma'$ which identifies the collars ($\iota_{\Sigma'} \circ
  T|_{N_{\Sigma}} = \iota_{\Sigma'}$) and such that the maps $g' \circ
  T$ and $g$ are thin homotopic.

  Composition of two morphisms $g$ and $g'$ is defined by gluing the
  outputs of $g$ to the inputs of $g'$ and is denoted by $g'\circ g$.
  This composition is associative because of the identification of
  diffeomorphic $X$-surfaces.  If $\gamma\colon S_n \to X$ is an
  object then the identity morphism is the $X$-surface
  $\Id_\gamma\colon C_n\to X$, recalling that $C_n$ is $n$ cylinders,
  given by composing the projection $C_n\to S_n$ with $\gamma$,
  because gluing $\Id_\gamma$ to an $X$-surface is thin-homotopic to
  the original $X$-surface.

  Disjoint union $\sqcup$ of $X$-surfaces makes $\thsx$ into a strict
  symmetric monoidal category (see Appendix \ref{ap:defns} for the
  definition of a symmetric monoidal category).  The unit for this
  monoidal structure is the empty $X$-surface. For objects
  $\gamma\colon S_n \ra X$ and $\gamma^\prime\colon S_{n^\prime} \ra
  X$ the symmetry structure isomorphism $\kappa\colon \gamma\sqcup
  \gamma^\prime \ra \gamma^\prime \sqcup \gamma$ is given by the
  cylinder $C_{n+ n^\prime} \ra S_{n+ n^\prime} \ra X$ where
  the first map is the projection and the boundary identification
  applies the appropriate permutation of boundary circles.
\end{defn}
Note that if we did not include thin-homotopy in the definition
then we would not have a category as there would not be any
identity morphisms.

Now we introduce the main definition of this section.

\begin{defn}
  A {\em rank one, smooth, thin-invariant field theory} for a smooth
  manifold $X$ is a symmetric monoidal functor $\ehqft$ (see
  Appendix~\ref{ap:defns}) from the thin homotopy surface category of
  $X$ to the category of one-dimension complex vector spaces with
  tensor product, satisfying the following smoothness condition.  If
  $g$ is a closed surface, then write it as $\lan g\ran$ to emphasise
  the fact that it is closed. Such a closed $X$-surface is an
  endomorphism of the empty object so $E\lan g\ran$ is a linear map on
  $\bC$ so can be identified with a complex number, this number is the
  {\em holonomy} of $\lan g\ran$ and will also be written $E\lan
  g\ran$.  The smoothness condition is then that there exists a closed
  3-form $c$ on $X$ such that if $v\colon V\to X$ is an
  $X$-three-manifold then \[ E\lan \partial v\ran = \exp(2\pi
  i{\textstyle\int_V v^*c}).  \] Two thin-invariant field theories are
  {\em isomorphic} if there is a monoidal natural isomorphism between
  them. If the three-form $c$ is zero we say that the thin-invariant
  field theory is\/ {\em flat}.
\end{defn}

It is possible to define higher rank thin-invariant field theories,
and these should be related to non-Abelian gerbes, but we will not
discuss them here.  For the rest of this paper ``thin-invariant field
theory'' will mean ``rank one, smooth thin-invariant field theory''.

Note that according the definition of a symmetric monoidal functor (see
Appendix~\ref{ap:defns}) a thin-invariant field theory $\ehqft$ comes
equipped with natural isomorphisms $\Phi_{\gamma,\gamma^\prime}\colon
E(\gamma)\otimes E(\gamma^\prime) \rightarrow E(\gamma \sqcup
\gamma^\prime)$ for each pair of objects $\gamma$ and
$\gamma^\prime$. These are symmetric, that is,
\[
\Phi_{\gamma^\prime,\gamma} \circ T = E(\kappa) \circ
\Phi_{\gamma,\gamma^\prime}
\]
where $T$ is the flip in $\vect_1$ and $\kappa\colon \gamma\sqcup
  \gamma^\prime \ra \gamma^\prime \sqcup \gamma$ is the symmetry
structure isomophism in the thin-homotopy surface category.

The definition of isomorphism of thin-invariant field theories
requires a natural transformation $\Psi\colon E\ra E^\prime$ such that
for each object $\gamma$, the map $\Psi_\gamma\colon E(\gamma) \ra
E^\prime (\gamma)$ is an isomorphism and for each pair of objects
$\gamma$ and $\gamma^\prime$
\[
\Psi_{\gamma \sqcup \gamma^\prime} \circ \Phi^E_{\gamma,\gamma^\prime}
= \Phi^{E^\prime}_{\gamma,\gamma^\prime} \circ
( \Psi_\gamma \otimes \Psi_{\gamma^\prime}).
\]

There is a\/ {\em trivial\/} thin-invariant field theory defined by
setting $E(\gamma) = \bC$ for all objects $\gamma$, setting $E(g):=\Id$
for all morphisms $g$, and taking $\Phi_{\gamma,\gamma^\prime}\colon
\bC \otimes \bC \rightarrow \bC$ to be the canonical identification.

\begin{prop}
  The set of isomorphism classes of thin-invariant field theories on a
  manifold $X$ form a group which will be denoted by $\Tift(X)$.
  Furthermore the flat thin-invariant field theories on $X$ form a
  subgroup.
\end{prop}
\begin{proof}
 Given thin-invariant field
  theories $E$ and $F$ there is thin-invariant field theory $E\otimes
  F$ formed by defining $(E\otimes F)(\gamma) := E(\gamma) \otimes
  F(\gamma)$ for objects, $(E\otimes F)(g) = E(g) \otimes F(g)$ for a
  morphism $g$ and $\Phi_{E\otimes F} = (\Phi_E\otimes \Phi_F) \circ
  T$ where $T$ is the flip. The three form $c_{E\otimes F}$ is equal
  to $c_E+c_F$.   The identity of this
  group is the trivial thin-invariant field theory.  The inverse
  $E^{-1}$ of $E$ is defined by setting $E^{-1}(\gamma)=(E(\gamma))^*=  \Hom (E(\gamma),\bC)$ for objects, $E^{-1}(g) = E(\op{g})^*$ for a
  morphism $g$ and $\Phi_{E^{-1}} = (\Phi^{-1}_E)^*$.
 \end{proof}

The next lemma is a useful property coming from the fact that we are
only considering the rank one case.
\begin{lem}\label{lem:closed}
  Suppose that $E$ is a thin-invariant field theory on the smooth
  manifold $X$.  If $g\colon\Sigma\to X$ is an endomorphism of the
  object $\gamma$ of $\thsx$ and $\lan g\ran$ is the closed
  $X$-surface obtained by identifying the inputs and outputs of $g$
  then
  \[ E(g)=E\lan g\ran\Id_\gamma.  \]
\end{lem}
\begin{proof}
  This is a standard argument in topological field theory.  The
  cylinder, thought as a cobordism from $\gamma\sqcup\gamma$ to the
  empty map, gives rise to a non-degenerate inner-product on
  $E(\gamma)$.  Evaluating $E\lan g\ran$ is the same as calculating
  the trace of $E(g)$ using this inner product.  The result follows
  from this because $E(g)$ is an endomorphism of a one-dimensional
  space.
\end{proof}

The following theorem gives a fundamental property of thin-invariant
field theories.
\begin{thm}\label{thm:cobinv}
A thin-invariant field theory is invariant under thin cobordism of morphisms.
\end{thm}
\begin{proof}
  Suppose that $E$ is a thin-invariant field theory.
  It suffices to show that if $g$ and $g'$ are thin-cobordant then
  $E(g')= E(\op g)^{-1}$, because, as $g$ is thin-cobordant to itself,
  we also get $E(g)= E(\op g)^{-1}$ and hence $E(g')= E(g)$.

  So suppose that $w$ is a thin-cobordism with $\partial w\cong
  \lan g'\circ \op{g}\ran$.  Then if $c$ is the three-form of $E$ we get
  that $\int_W w^*c=0$ as $dw$ everywhere has rank two, so $E\lan g'
  \circ \op{g}\ran=1$.  By using the previous lemma we find $E(g'
  \circ \op{g})=\Id$ from whence  $E(g')\circ E(\op{g})=\Id$, and
  $E(g')= E(\op g)^{-1}$ as required.
\end{proof}
This means that a thin-invariant field theory descends to a symmetric monoidal
functor on the thin-cobordism surface category of $X$, the category
obtained by replacing ``thin-homotopic'' by ``thin-cobordant'' in the
above definition. One fundamental property of this category is that it
is a groupoid, unlike the thin-homotopy category. This is proved by
the proposition in Section~2 and we get the important relation $E(\op
g)=E(g)^{-1}$.


\section{On flat thin-invariant field theories and homotopy quantum field
theories} We will elucidate the connection between thin-invariant
field theories and homotopy quantum field theories by showing that
a certain subset of homotopy quantum field theories, the rank one,
normalised ones, are the same as flat thin-invariant field
theories.

Recall that the homotopy surface category is defined by replacing
the term ``thin-homotopic'' by the term ``homotopic'' in the
definition of the thin-homo\-topy surface category.  A
1+1-dimensional homotopy quantum field theory on a space is a
symmetric monoidal functor from the homotopy surface category of
the space to the category of vector spaces.  This is a slight
variation on Turaev's original definition better suited to spaces
with a non-trivial second homotopy group. We will be interested in
rank one homotopy quantum field theories, that is those which are
functors to the subcategory of one dimensional vector spaces. In
what follows, HQFT means rank one, 1+1-dimensional homotopy
quantum field theory.

A flat thin-invariant field theory is one whose three-form is zero, so
it descends to a functor on the homotopy surface category and gives
rise to an HQFT.  However not all HQFTs arise in this way, as is
illustrated by the case of a point.  An HQFT on a point is the same thing as a topological
quantum field theory and a rank one topological quantum field theories
is determined by the invariant $\omega\in\bC$ of the two-sphere, the
genus $l$ surface having invariant $\omega^{2-l}$.  On the other hand,
all thin-invariant field theories on a point are trivial, as all
surfaces are cobordant.  We wish to compensate for this, so we make
the following definition.

\begin{defn} An HQFT on a space $X$ together with a point in $X$ gives rise to
a topological quantum field theory by considering the constant
$X$-surfaces to the point.  Note that if two points are in the
same connected component then the  topological quantum field
theories induced in this way are isomorphic. An HQFT is {\em
normalised} if for every point in $X$ the induced topological
quantum field theory is  trivial.
\end{defn}

The key property of a
normalised HQFT is that holonomy of a closed $X$-surface only depends
on the homology class of the $X$-surface.  This is the content of the
following proposition.
\begin{prop}\label{Prop:HomologyHolonomy}
For a rank one, normalised, 1+1-dimensional homotopy quantum field
theory, the holonomy of a closed $X$-surface $g\colon\Sigma\to X$
depends only on the
homology class $g_*[\Sigma]\in H_2(X;\bZ)$.
\end{prop}
We delay the proof until after the next theorem.

\begin{thm}\label{thm:flatequalsnormalixed}
  Every rank one normalised 1+1-dimensional homotopy
  quantum field theory can be considered as a flat
  thin-invariant field theory and\/ {\it vice versa}.
\end{thm}
\begin{proof}
  The discussion earlier shows that every flat thin-invariant field
  theory can be thought of as a normalised HQFT and we now prove that
  the converse also holds.

  To show that every normalised HQFT comes from a thin-invariant field
  theory we need to show that it satisfies the three-form condition
  with the three-form equal to zero, ie.\ if $v\colon V\to X$ is an
  $X$-three-manifold then
  $H(\partial v)=1$, but that is true because $H(\partial v)$ only
  depends on the homology class of $\partial v$, which is zero as it
  is cobordant to the empty manifold.
\end{proof}

The remainder of this section is dedicated to the proof of
Proposition~\ref{Prop:HomologyHolonomy}. First we need some results
about local surgery.

\begin{defn}
  An $X$-surface $g'\colon \Sigma' \to X$ is said to be obtained\/
  {\em by a local surgery\/} from $g\colon \Sigma \to X$ if $\Sigma'$
  with two discs removed is diffeomorphic to $\Sigma$ with a cylinder
  removed, the maps $g$ and  $g'$ agree on the diffeomorphic parts and
  they map the discs and cylinder mapped to a single point of $X$.
 (See Figure~1.)
\end{defn}
\begin{figure}
\cl{\vpic{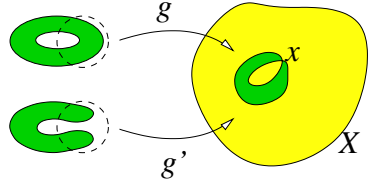}{1.9in}}
\caption{An example of a local surgery.  Both parts contained in the
  dashed circles are mapped to the point $x$.}
\end{figure}

 Informally this says that while the surfaces may be topologically
 different, they differ only at the inverse image of a point in $X$.

\begin{lem}
  \label{Lemma:BordantIffSequence}
  Two closed $X$-surfaces are cobordant if and only if they can be
  connected by a sequence of homotopies, local surgeries, and disjoint
  unions with null-homotopic spheres.
\end{lem}
\begin{proof}
  Suppose that $w\colon W\to X$ is a cobordism from $g\colon \Sigma
  \to X$ to $g'\colon \Sigma' \to X$.  By standard results in Morse
  theory \cite{Milnor:hCobordismBook}, we can pick a Morse function
  $f\colon W\to [0,1]$ such that $f^{-1}(0)=\Sigma$ and
  $f^{-1}(1)=\Sigma'$.  Let $0<f_1<\dots<f_c<1$ be the critical values
  of $f$, then the restricted maps $w\colon f^{-1}\left( (f_i,f_{i+1})
  \right) \to X$ for $i=1\dots c-1$ give a sequence of homotopies.
  The critical points of index zero and three correspond to the
  addition and deletion of null-homotopic two-spheres, while those of
  index one and two correspond to local surgeries.

  Conversely, if we have a sequence of such alterations connecting
  $g$ and $g'$ then this gives rise to a bordism by reversing the
  above procedure.
\end{proof}
\begin{lem}
  \label{Lemma:UnchangedBySurgery}
  The holonomy of a rank one, normalised homotopy quantum field theory
  around a closed  $X$-surface is
  unchanged by local surgery and by the disjoint union with a
  null-homotopic two-sphere.
\end{lem}
\begin{proof}
  Firstly, if $s\colon S^2\to X$ is a null-homotopic map then it is
  homotopic to the constant map $S^2\to \{*\}$ for some point $*\in
  X$.  It follows that $H(S^2\to \{*\})=1$ as $H$ is normalised, ie.\
  induces the trivial topological quantum field theory.  Thus if $g$
  is any $X$-surface then
\[H(g\sqcup s)=H(g)H(s)=H(g)H(S^2\to \{*\})=H(g).\]

Secondly, consider the union of two discs, $D\sqcup D$ and the
cylinder $C$, as surfaces with two inputs and no output.  As the
induced TQFT of $H$ is trivial we have $H(D\sqcup D\to \{*\})=H(C\to
\{*\})$ as maps $H(S^1\sqcup S^1)\to \bC$.

  If $g'\colon \Sigma' \to X$ is obtained from $g\colon \Sigma \to X$
  by a local surgery then let $\Omega$ be the surface with two
  outgoing boundary component such that $\Omega \cup \op D \cup \op D
  \cong \Sigma$ and $\Omega \cup C\cong \Sigma'$ with
  $g|_\Omega=g'|_\Omega$ and such that $g|_{D\sqcup D}$ and $g|_C$ are constant
  maps to the point $\{*\}$.  Thus
  \begin{align*}
    H\lan g\ran&=H\left\lan (D\sqcup D\to \{*\})\circ g|_\Omega \right\ran
        =H\left( D\sqcup D\to \{*\} \right)\circ H\left(g|_\Omega\right)\\
          &=H\left( C\to \{*\} \right)\circ H\left(g|_\Omega\right)
        =H\left\langle( C\to \{*\})    \circ g'|_\Omega \right\rangle\\
        &=H\lan g'\ran.
  \end{align*}
Which is what was required.
\end{proof}

We can now prove Proposition~\ref{Prop:HomologyHolonomy}, which stated
that the holonomy of a normalised HQFT depends only on the
homology of $X$-surfaces.
\begin{proof}[Proof of Proposition~\ref{Prop:HomologyHolonomy}]
  Suppose that  $g\colon \Sigma \to X$ and $g'\colon \Sigma' \to X$ are
  homologous closed $X$-surfaces, in the sense that
  $g_*[\Sigma]=g'_*[\Sigma']\in H_2(X,\bZ)$.  We need to show that $H\langle
  g\rangle =H \langle g' \rangle$.

  By Lemma \ref{Lemma:SecondHomologyEqualsBordism}, $g$ is bordant to
  $g'$, so Lemma~\ref{Lemma:BordantIffSequence} implies that there is
  a sequence of homotopies, local surgeries, and disjoint unions with
  two-spheres connecting $g$ to $g'$.  The definition of an HQFT
  ensures that the holonomy does not change under homotopy, and
  Lemma~\ref{Lemma:UnchangedBySurgery} ensures that it is unchanged
  under the latter two as well.  Thus $H(g)=H(g')$ as required.
\end{proof}


\section{Examples of thin-invariant field theories}\label{sec:examples}
In this section we present a number of examples of thin-invariant
field theories.

\subsection{Manifolds with trivial first homology}\label{ex:scm}
The first example applies to spaces with trivial first integral
homology group.  For such a space we build a thin-invariant field
theory starting from a gerbe with connection. The construction was
partially inspired by \cite{Freed:QuantumGroups}.

Let $X$ be a smooth finite dimensional manifold with $H_1(X,\bZ)$
trivial, and let $S$ be the holonomy of a gerbe with connection.  If
$\gamma$ is an object in the thin-cobordism category $\thsx$ then
Lemma \ref{Lemma:SecondHomologyEqualsBordism} implies that $\gamma$ is
null cobordant. Define the one dimensional vector space $E(\gamma)$ to
be the space of complex linear combinations of null cobordisms of
$\gamma$ modulo a relation involving the gerbe holonomy $S$:
\[
  E(\gamma):=\bC\Hom_{\thsx}(\emptyset,\gamma)\bigm / \{h_1=S\langle
  h_1\circ \op{h_2}\rangle h_2\}.
\]
This is clearly one-dimensional.

If $g$ is a morphism in $\thsx$ from $\gamma$ to $\gamma'$
then define
$  E(g)\colon E(\gamma)\to E(\gamma')$ by $E(g)h:= g\circ h$.
This is well defined on $E(\gamma)$ because of the following:
\begin{align*}
E(g)(S\lan h_1\circ\op{h_2}\ran h_2)
  &= S\lan h_1\circ\op{h_2}\ran g\circ h_2
   = S\lan g\circ h_1\circ\op{h_2}\circ \op{g}\ran g\circ h_2\\
  &= g\circ h_1 = E(g)h_1.
\end{align*}
Functoriality is immediate: $E(g\circ g')=E(g)\circ E(g')$.

To show thin-invariance we need to show that if $g$ is thin-cobordant
to $g'$ then $E(g)=E(g')$.  If $g$ is thin-cobordant to $g'$ then by
Proposition~\ref{prop:holonomythininvariance}, $S\lan\op{g'}\circ
g\ran=1$ and we find
\begin{align*}
  E(g)h&=g\circ h = S\lan g\circ h\circ\op{h}\circ \op{g'}\ran g'\circ h
      = S\lan h\circ\op{h}\circ \op{g'}\circ g\ran g'\circ h\\
    & = S\lan h\circ \op{h}\ran S\lan \op{g'}\circ g\ran  g'\circ h
      = g'\circ h
      = E(g')h.
\end{align*}
It follows from this invariance that $E$ also respects the identity
maps.

The smoothness condition is automatically satisfied by the curvature
three-form of the gerbe.

To show that $E$ is symmetric monoidal it is necessary to show that
there are symmetric natural
isomorphisms:
\[
  \Phi_{\gamma,\gamma'}\colon E(\gamma)\otimes E(\gamma')\stackrel{\cong}{\longrightarrow}
     E(\gamma \sqcup \gamma').
\]
Define $\Phi_{\gamma,\gamma'}(h\otimes h^\prime) = h\sqcup h^\prime $ and
note this
is well defined since
\begin{align*}
\Phi_{\gamma,  \gamma'}(S\lan h_1\circ\op{h_2}\ran h_2\otimes S\lan
h^\prime_1\circ\op{h^\prime_2}\ran h^\prime_2) &S\lan h_1\circ\op{h_2}\ran S\lan h^\prime_1\circ\op{h^\prime_2}\ran
h_2\sqcup h^\prime_2 \\
& = S\lan (h_1\circ\op{h_2})\sqcup ( h^\prime_1\circ\op{h^\prime_2})\ran
h_2\sqcup h^\prime_2 \\
&= S\lan (h_1\sqcup h^\prime_1) \circ (\op{h_2\sqcup
h^\prime_2})\ran h_2\sqcup h^\prime_2 \\
&= h_1\sqcup h^\prime_1 .
\end{align*}
Moreover, if $g\in \Hom_{\thsx} (\gamma_1, \gamma_2)$ and $g'\in
\Hom_{\thsx} (\gamma^\prime_1, \gamma^\prime_2)$ then
\begin{align*}
E(g\sqcup g') \circ \Phi_{\gamma_1, \gamma^\prime_1}(h\otimes h^\prime) &E(g\sqcup g') (h\sqcup h^\prime) = (g\sqcup g')\circ (h\sqcup
h^\prime)\\
& = (g\circ h)\sqcup (g'\circ h^\prime) \Phi_{\gamma_2,\gamma^\prime_2}(g\circ h\otimes g'\circ h^\prime) \\
& = \Phi_{\gamma_2,\gamma^\prime_2}(E(g) h\otimes E(g')h^\prime) \\
&=\Phi_{\gamma_2,\gamma^\prime_2} \circ (E(g)\otimes E(g'))(h\otimes h^\prime),
\end{align*}
proving that the $\Phi_{\gamma,\gamma^\prime}$ are natural.  Let $T$
be the flip $E(\gamma) \otimes E(\gamma^\prime) \rightarrow
E(\gamma^\prime) \otimes E(\gamma)$ and $\kappa\in\Hom_{\thsx}
(\gamma\sqcup\gamma^\prime, \gamma^\prime\sqcup\gamma)$ be the
symmetric structure isomorphism for $\gamma$ and $\gamma^\prime$. Then
\begin{align*}
\Phi_{\gamma^\prime,\gamma} \circ T (h\otimes h^\prime)&=\Phi_{\gamma^\prime,\gamma}(h^\prime\otimes  h) = h^\prime \sqcup h
 = S\lan (h \sqcup h^\prime)\circ  \op{\kappa \circ (h \sqcup h^\prime)}\ran
   \kappa\circ (h \sqcup h^\prime) \\
 &  = E(\kappa)(h \sqcup h^\prime)
= E(\kappa) \circ \Phi_{\gamma,\gamma^\prime}(h\otimes  h^\prime).
\end{align*}

 This proves that $E$ is a thin-invariant field theory.

\subsection{Gerbes with trivial Dixmier-Douady class}\label{ex:ddz}
The second example does not require any restrictions on the manifold
$X$, and builds a thin-invariant field theory from a gerbe with
connection whose Dixmier-Douady class is zero.  By the exact sequence
of Section~\ref{Section:GerbeHolonomy} such a gerbe may be
represented, non-uniquely, by a 2-form $\omega$, with the holonomy
around a closed $X$-surface $g\colon\Sigma\ra X$ given by $\exp (2\pi
i \int_\Sigma g^*\omega)$.  Now define a thin-invariant field theory
by setting $E(\gamma) = \bC$ for each object $\gamma$ in $\thsx$, and
for each morphism $g\colon\Sigma\ra X$ defining $E(g)\colon \bC \ra
\bC$ to be multiplication by $\exp (2\pi i \int_\Sigma g^*\omega)$.
The monoidal structure $\Phi_{\gamma,\gamma^\prime}
\colon\bC\otimes\bC\ra\bC$ is the canonical isomorphism.
Thin-invariance follows from Stokes' Theorem.  It is evident that the
holonomy of this thin-invariant field theory is the same as that of
the original gerbe.

\subsection{Flat thin-invariant field theories from two-cocycles}
\label{ex:turaev}
The third example uses the identification of flat thin-invariant field
theories with normalised HQFTs
(Theorem~\ref{thm:flatequalsnormalixed}) to get examples of flat thin-invariant
field theories from a construction of Turaev
\cite{Turaev:HomotopyFieldTheoryInDimension2}.  Note that we use a
slightly different convention to Turaev to ensure that we get the
correct holonomy, and not its inverse.  Let $\theta\in
C^2(X,\bC^\times)$ be a two-cocycle, and define $E^\theta (\gamma)$ for an
object $\gamma\colon S_m\to X$ by taking all one-cycles which
represent the fundamental class of $S_m$ and
quotienting by a certain relation:
\[
  E^\theta (\gamma):=\bC\biggl\{ a\in C_1(S_m)\biggm |\!
  [a]=[S_m]\biggr\} \biggm /\!
  \bC\biggl\{ a-\gamma^*\theta(e)b\biggm|\!\!
   \begin{array}{l} e\in C_2(S_m)\\ \partial e=a-b
   \end{array}
  \biggr\}.
\]
Write $|a|$ for the equivalence class in $E^\theta (\gamma)$ of the
one-cycle $a$.  To a cobordism $g\colon \Sigma\to X$ from $\gamma_0$
to $\gamma_1$ we need to associate a linear map.  This is done by
picking a singular two-cycle representative $f\in C_2(\Sigma)$ of the
fundamental class $[\Sigma]\in H_2(\Sigma,\partial \Sigma)$.  Then
$E^\theta (g)\colon E^\theta (\gamma_0) \to E^\theta ({\gamma_1})$ is
defined by $E^\theta (g)|a_0|:= g^*\theta(f) |a_1|$, where $\partial f
= a_0 - a_1$.  For objects $\gamma_0$ and $\gamma_1$ the monoidal
structure map $\Phi_{\gamma_0,\gamma_1}\colon E^\theta
(\gamma_0)\otimes E^\theta (\gamma_1)\rightarrow E^\theta
(\gamma_0\sqcup \gamma_1)$ is defined such that
$\Phi_{\gamma_0,\gamma_1}(|a_0|\otimes |a_1|):=|a_0\sqcup a_1|$.

Turaev shows that this is well-defined and gives a normalised HQFT,
hence, by Theorem~\ref{thm:flatequalsnormalixed}, we get a flat
thin-invariant field theory.  In actual fact this gives rise to a
group homomorphism from $H^2(X,\bC^\times)$ to the group of flat
thin-invariant field theories up to equivalence, as if two
two-cocycles differ by a coboundary then the
thin-invariant field theories are non-canonically isomorphic in the
following manner.  If $\theta=\theta'+\delta f$ where $f\in C^1(X)$
then for an object $\gamma$ in $\thsx$ define $\Psi_\gamma\colon
E^\theta (\gamma) \ra E^{\theta^\prime}(\gamma)$ by
$\Psi_\gamma(|a|):=(f(\gamma_*a))^{-1}|a|$: it transpires that $\Psi$
is a natural transformation giving an isomorphism of thin-invariant
field theories.  We have a group homomorphism because the
theory $E^{\theta_1 + \theta_2}$ constructed from $\theta_1 +
\theta_2\in C^2(X)$ is isomorphic to $E^{\theta_1}\otimes E^{\theta_2}$.


\section{Thin-invariant field theories and gerbes}
\label{Section:TIFTsAndGerbes}
The goal of this section is to show that a gerbe with a connection is
the same thing as a thin-invariant field theory.


\subsection{Ext groups and monoidal functors}
To prove the main theorem we are going to need an aside on $\Ext$
groups.  We will start with a little reminder.  If $\Gamma$ and $A$
are abelian group then $\Ext(\Gamma, A)$ is the set of all abelian
extensions of $\Gamma$ by $A$, that is to say, all abelian groups
$\hat\Gamma$ which fit into an exact sequence $0\to A\to \hat\Gamma\to
\Gamma\to 0,$ modulo some suitable notion of equivalence.  Similarly
the group cohomology group $H^2_{\text{gp}}(\Gamma,A)$ can be
identified with the set of central extensions of $\Gamma$ by $A$, that
is those $\hat\Gamma$ as above which are not necessarily abelian, but
in which $A$ is a central subgroup.

We will be interested in the case that $A$ is $K^\times$, the group of
units of a commutative ring $K$.  In this case a $K^\times$-extension
of $\Gamma$ is like a $K$-line bundle over the discrete space
$\Gamma$.  We will need the notion of a $K^\times$-torsor which is
just another name for a principal $K^\times$-homogeneous space, ie.\ a
space with a transitive and free $K^\times$-action.  Of course such a
space is homeomorphic to $K^\times$, but in general there will be no
canonical homeomorphism.  The collection of $K^\times$-torsors forms a
symmetric monoidal category, in which the morphisms are the maps
commuting with the action, the monoidal product is given by $A\cdot
B:=A\times_{K^\times}B$, and the unit object is just $K^\times$.

\begin{lem}
   \label{Lemma:KTorsorsKModules}
  If $K$ is a commutative ring then there is an equivalence of
  symmetric monoidal categories between the category of
  $K^\times$-torsors and the category of rank one $K$-modules.
\end{lem}
\begin{proof}
  The functor $\Ktor\to\Kmod$ is given by $A\mapsto K\times_{K^\times}A$,
  and the functor $\Kmod \to \Ktor$ is given by $M\mapsto \{
  \text{generators of $M$ over $K$}\}$.  The reader is left to fill in
  the uninspiring details.
\end{proof}

Given an abelian group $\Gamma$ we can form a strict symmetric monoidal
category $\underline\Gamma$ whose objects are the elements of the
group, whose morphisms are just the identity morphisms and whose
monoidal structure is just given by the group multiplication.  In the
next section we will be using $\underline{H_1(X,\bZ)}$, which is just
the category consisting of the ``connected components'' of the thin-cobordism
category. Now we will relate functors from $\underline \Gamma$ with
extensions of the group $\Gamma$.

\begin{lem}
\label{lem:ext}
  Suppose that $\Gamma$ is an abelian group and $K$ is a commutative
  ring.  There is a naturual bijection between
  $\Ext(\Gamma,K^\times)$, the abelian extensions of $\Gamma$ by
  $K^\times$, and the set of symmetric monoidal functors from the
  category $\underline\Gamma$ defined above to the category of rank
  one $K$-modules.
\end{lem}
\begin{proof}
  Firstly, in view of the previous lemma, we can equivalently consider
  symmetric monoidal functors from $\underline\Gamma$ to the category
  of $K^\times$-torsors.

  So suppose that $E\colon \underline\Gamma\to \Ktor$ is a symmetric
  monoidal functor.  Define $\hat\Gamma:=\bigcup_{x\in\Gamma}E(x)$.
  We need to show that this is an abelian extension of $\Gamma$ by
  $K^\times$.  The multiplication comes from the monoidal structure
  $E(x)\cdot E(y)\to E(xy)$, this is associative because of the
  associativity axiom for monoidal structure.  Inverses exist for the
  following reason: if $\alpha\in \hat\Gamma$ lives in $E(x)$ then
  there is a map $\Phi_x\colon E(x)\cdot E(x^{-1})\to K^{\times}$,
  pick any element $\beta\in E(x^{-1})$ and take $\alpha^{-1}$ to be
  $(\Phi_x(\alpha,\beta))^{-1}\beta$.  There is the obvious quotient
  map $p\colon \hat\Gamma\to\Gamma$, which is automatically a group
  homomorphism, and there is the inclusion homomorphism
    $K^\times{\cong}E(1)\hookrightarrow\hat\Gamma$
  coming from the unit axiom for a symmetric monoidal
  functor.  The symmetric axiom then gives that $\hat\Gamma$ is abelian.

  Conversely, if $0\to K^\times\stackrel{i}\to \hat\Gamma\stackrel{p}
  \to \Gamma\to 0$ is an abelian extension of $\Gamma$ by $K^\times$
  then define $E\colon\underline \Gamma \to \Ktor$ by
  $E(x):=p^{-1}(x)$.  Note that $K^\times$ acts transitively and
  freely on $p^{-1}(x)$ via $i$.  The structure maps $E(x)\cdot
  E(y)\to E(xy)$ come directly from the product in $\hat\Gamma$ and
  these are symmetric because $\hat\Gamma$ is abelian.  Isomorphic
  extensions give isomorphic functors.
\end{proof}
Note that if we replace the phrase ``abelian extension'' by ``central
extension'' and ``symmetric monoidal functor'' by ``monoidal functor''
in the above proof then we get a bijection between
$H^2_{\text{gp}}(\Gamma,K^\times)$ and monoidal functors from $\underline \Gamma$
to the category of rank one $K$--modules.


\subsection{The Main Theorem}
Now we can prove the main result of this paper.

\begin{thm}
\label{thm:GerbeEqualsTIFT}
    \setcounter{test}{\value{thm}}\addtocounter{test}{-1}
    \refstepcounter{test}\label{thmnumber:GerbeEqualsTIFT}
  On a smooth finite dimensional manifold, there is
  an isomorphism from the group of thin-invariant field theories (up to
  equivalence) to the group of gerbes with connection (up to equivalence).
\end{thm}

\begin{proof}
We will show that the holonomy $S\colon \Tift (X) \rightarrow
\widehat H^3(X)$ is an isomorphism.

Firstly to show that $S$ is injective we will identify the kernel of
$S$ with the set of symmetric monoidal functors
$\underline{H_1(X,\bZ)} \rightarrow \text{Vect}_{1}$. Then we can use
Lemma~\ref{lem:ext} to identify this set with $\Ext
(H_1(X,\bZ),\bC^\times)$ which we know to be trivial as $\bC^\times$
is a divisible group.

To identify the kernel of $S$ with the collection of symmetric
monoidal functors $\underline{H_1(X,\bZ)} \rightarrow
\text{Vect}_{1}$ we proceed as follows.  Suppose $E\in \Ker(S)$,
then we will construct a symmetric monoidal functor
$\cE\colon\underline{H_1(X,\bZ)} \rightarrow \text{Vect}_{1}$.
Suppose $\gamma_1$ and $\gamma_2$ are objects in $\thsx$, and
suppose $g,g^\prime\in \Hom_{\thsx} (\gamma_1, \gamma_2)$.  Since
$E$ has trivial holonomy, we have $E\lan g\circ \op{g^\prime} \ran
= 1$ and it follows from Lemma \ref{lem:closed} that $E(g)\circ
E(\op{g^\prime}) = \Id$ and hence from the discussion after
Theorem~\ref{thm:cobinv} that $E(g^\prime) = E(g)$. Thus, for each
cobordant $\gamma_1$ and $\gamma_2$ there is a canonical
identification of $E(\gamma_1)$ with $E(\gamma_2)$. By Lemma
\ref{Lemma:SecondHomologyEqualsBordism}, $\gamma_1$ and $\gamma_2$
are cobordant if and only if they belong to the same homology
class, and we can therefore associate in a natural way a one
dimensional vector space $\cE(x)$ to each homology class $x\in
H_1(X,\bZ)$.

There are natural isomorphisms $\Phi_{x,x^\prime}\colon \cE(x)\otimes
\cE(x^\prime) \rightarrow \cE(x+x^\prime)$ obtained by choosing $\gamma$
to represent $x$ and $\gamma^\prime$ to represent $x^\prime$, so that
$\cE(x)\cong E(\gamma)$ and $\cE(x^\prime)\cong E(\gamma^\prime)$, and
then setting $\Phi_{x,x^\prime} = \Phi_{\gamma,\gamma^\prime}$.  It
follows from properties of $\Phi_{\gamma,\gamma^\prime}$ that these
are well defined natural isomorphisms.  Thus, $\cE$ is a
symmetric monoidal functor.

If the functor just defined is isomorphic to the trivial one then
the field theory giving rise to it must also be trivial.
Conversely, any monoidal functor from $\underline{H_1(X,\bZ)}$ to
$\text{Vect}_{1}$ can be extended to a field theory with trivial
holonomy. Thus there is a bijection of $\Ker (S)$ with symmetric
monoidal functors $\underline{H_1(X,\bZ)}\rightarrow
\text{Vect}_{1}$ as required. As explained at the beginning of the
proof, this shows that $\Ker (S)=0$.

The second step is to show that $S$ is surjective.
Example~\ref{ex:scm} shows that if $X$ is simply connected then every
gerbe is the image under $S$ of some thin-invariant field theory.
Similarly Example~\ref{ex:ddz} shows that for arbitrary $X$, every
gerbe with zero Dixmier-Douady class comes, via $S$, from a
thin-invariant field theory.

For the general case suppose we have a gerbe holonomy $\cS$.
Let $M$ be a smooth
manifold which is dim$(X)+1$-homotopy equivalent to $K(\bZ,3)$. By
identifying $H^3(X,\bZ)$ with homotopy classes of maps $X\rightarrow
M$ we can choose a smooth map $f\colon X\rightarrow M$ representing
the Dixmier-Douady class of $\cS$.
Now let $\cS_1$ be a gerbe over $M$ whose Dixmier-Douady class
is the generator of $H^3(M,\bZ)$. Since $M$ is simply connected we can apply
Example~\ref{ex:scm} to obtain a thin-invariant field theory $E_1$
over $M$ such that $S_{E_1}=\cS_1$.  The gerbe-holonomy $\cS\otimes
f^*\cS_1^{-1}$ has Dixmier-Douady class zero and we can apply Example
\ref{ex:ddz} to obtain a thin-invariant field theory $E_0$ over $X$.
Finally, the thin-invariant field theory $E=E_0\otimes f^*E_1$
satisfies $S_E = \cS$.
\end{proof}

Note that the second half of this proof could have been done more
neatly if we could classify the gerbe holonomy by a smooth map
$X\to BBS^1$, but this does not seem to be in
\cite{Gajer:GeometryOfDeligneCohomology}.

Combining this with Theorem~\ref{thm:flatequalsnormalixed}, we obtain
the identification of normalised HQFTs with flat gerbes:
\begin{thm}
    \setcounter{test}{\value{thm}}\addtocounter{test}{-1}
    \refstepcounter{test}\label{thmnumber:FlatGerbeEqualsNHQFT}
  On a finite dimensional manifold, the group of normalised rank one homotopy
  quantum field theories (up to equivalence) is isomorphic to the
  group of flat gerbes (up to equivalence).
\end{thm}

The definition of thin-invariant field theory requires the functor
$\ehqft$ to be symmetric. In view of the motivation of this paper,
namely to reconcile homotopy quantum field theories and gerbes, this
is an entirely natural assumption to make. However, one can drop this
assumption, to get\/ {\em non-symmetric\/} thin-invariant field
theories.  In this case there is an analogue of Theorem~\ref{thm:GerbeEqualsTIFT}.
\begin{thm}
There is a split short exact sequence
\[
0 \rightarrow H^2_{\text{gp}}(H_1(X,\bZ), \bC^\times) \rightarrow
\left\{ \begin{array}{c}\mbox{non-symmetric}\\
                           \mbox{thin-invariant}\\
                           \mbox{field theories on $X$}
                         \end{array}\right\}
    \rightarrow  \widehat H^3(X) \rightarrow 0.
\]
\end{thm}
\begin{proof}
  The proof of Theorem~\ref{thm:GerbeEqualsTIFT} goes through almost
  exactly, by replacing the $\Ext$ group with the cohomology group, as
  in the comment after Lemma~\ref{lem:ext}. The splitting comes from
  the fact that we have already identified the group $\widehat H^3(X)$ with
  (symmetric) thin-invariant field theories.
\end{proof}


\section{Normalised homotopy quantum field theories and flat gerbes}
\label{Section:H2EqualsKHQFT}

Homotopy quantum field theories can be defined over rings other
than $\bC$, so let $K$ be any commutative ring, and recall that
HQFT is used to mean rank one 1+1-dimensional homotopy quantum
field theory.  Turaev's construction in Example~\ref{ex:turaev}
can be generalised to give a map $$ \tau\colon H^2(X, K^\times)
\rightarrow \{\text{normalised $K$-HQFTs}\}.$$ In
\cite{Turaev:HomotopyFieldTheoryInDimension2}, using his
classification of homotopy quantum field theories for
Eilenberg-MacLane spaces in terms of crossed algebras, Turaev
proved that when $X$ is a $K(\pi,1)$ and $K$ is a field of
characteristic zero then $\tau$ is an isomorphism.  We generalize
this in the following manner.

\begin{thm}
    \setcounter{test}{\value{thm}}\addtocounter{test}{-1}
    \refstepcounter{test}\label{thmnumber:H2EqualsKHQFT}
  Let $K$ be a commutative ring, and $X$ be a path
  connected topological space.  Then Turaev's construction gives an
  isomorphism between the group $H^2(X,K^\times)$ and the group of
  normalised, rank one homotopy quantum field theories defined over
  $K$.
\end{thm}

\begin{proof}
The proof will proceed like so.  We will construct the following diagram.
\[\xymatrix{
\Ext(H_1(X,\bZ),K^\times)\;\; \ar @{^{(}->}[r]^-{\iota} \ar@{=}[d] & H^2(X, K^\times)
\ar@{->>}[r] \ar[d]^\tau & \Hom (H_2(X,\bZ), K^\times)\ar@{=}[d] \\
\Ext(H_1(X,\bZ),K^\times)\;\; \ar @{^{(}->}[r] & \{\mbox{normalised
  HQFTs}\} \ar[r]^{S} &
\Hom (H_2(X,\bZ), K^\times)
}\]
We will show that it is commutative and then invoke the Five Lemma to
deduce that the map $\tau$ is an isomorphism.

First then we describe the morphisms in the diagram.  The top row is
the short exact sequence from the Universal  Coefficient
Theorem.  The map $\tau$ is Turaev's construction described in
Example~\ref{ex:turaev}.  The map $S$ is the holonomy map which is well
defined by  Proposition~\ref{Prop:HomologyHolonomy}.  The kernel of
$S$ is $\Ext(H_1(X,\bZ),K^\times)$ by the proof of
Theorem~\ref{thm:GerbeEqualsTIFT}.

Now we consider the commutativity of the diagram.  The right-hand
square is commutative by the definition of $\tau$.  The commutativity
of the left-hand square will take up the rest of the proof.

We need to describe the inclusion
$\iota\colon\Ext(H_1(X,\bZ),K^\times)\hookrightarrow H^2(X,K^\times)$,
in which we are considering $\Ext(H_1(X,\bZ),K^\times)$ as the group
of equivalence classes of extensions of abelian groups.  Let $B_1$ and
$Z_1$ be the groups of one-boundaries and one-cycles on $X$, then
$H_1(X,\bZ)$ has the following free-resolution: $ B_1\hookrightarrow
Z_1\stackrel{q}{\twoheadrightarrow} H_1(X,\bZ)$.  Now given $\epsilon$ an
abelian extension $ K^\times \hookrightarrow \hat{\Gamma}
\stackrel{p}{\twoheadrightarrow} H_1(X,\bZ)$, we can lift the morphism $q$ to
a morphism $\hat q\colon Z_1\rightarrow \hat{\Gamma}$.  Now we can define
$\theta\in C^2(X,K^\times)$ by
\[\theta_\epsilon(e):= \hat q(\partial e);
       \quad \text{for }e\in  C_2(X,K^\times).\]
The map $\theta$ is immediately seen to be a cocycle and it is
$K^\times$-valued because $\partial e\in B_1$ so $q(\partial e)=0$.
It transpires that the cohomology class $[\theta]$ is
precisely $\iota([\epsilon])\in H^2(X,K^\times)$.

Here we will take a slightly different but equivalent and more
convenient point of view of HQFTs, which is entirely analogous to
thinking of\/ {\em principal $\bC^\times$-bundles\/} rather than\/
{\em complex line bundles}.  We will think of HQFTs as associating
to an $X$-one-manifold a $K^\times$-torsor rather than a rank one
$K$-module.  In view of Lemma~\ref{Lemma:KTorsorsKModules} this
does not alter anything.  From this point of view, the HQFT
associated to $\theta$ by $\tau$, which will be denoted
$E^\theta$, is defined on $\gamma\colon S_m\to X$ by
\[
  E^\theta (\gamma):=K^\times\cdot\biggl\{ a\in C_1(S_m)\biggm |\!
  [a]=[S_m]\biggr\} \biggm /\!
  \biggl\{ a=\gamma^*\hat q(\partial e)b\biggm|\!\!
   \begin{array}{l} e\in C_2(S_m)\\ \partial e=a-b
   \end{array}
  \biggr\}.
\]
Remember that we use the notation $|a|$ for the equivalence class of
$a$ in $ E^\theta (\gamma)$.

The homomorphism sending the extension $\epsilon$ to $E^\theta$ gives the
composition of two of the sides of the left-hand
square in the diagram.  The composition of the other two is got by
looking at the proof of Theorem~\ref{thm:GerbeEqualsTIFT}, by which
we see that the HQFT, $E^\epsilon$, obtained from $\epsilon$,
is as follows.  Recalling that $\epsilon$ is the abelian extension $ K^\times
\hookrightarrow \hat{\Gamma} \stackrel{p}{\twoheadrightarrow}
H_1(X,\bZ)$, to each object $\gamma\colon S_m\rightarrow X$ we
associate the $K^\times$-torsor
$E^\epsilon(\gamma):=p^{-1}([\gamma])$, where $[\gamma]\in H_1(X,\bZ)$
denotes the class represented by $\gamma$.  The morphisms are mapped
to identities, and the symmetric monoidal structure
$\Phi^\epsilon_{\gamma,\gamma^\prime}\colon E^\epsilon(\gamma)\otimes
E^\epsilon (\gamma^\prime)\rightarrow
E^\epsilon(\gamma\sqcup\gamma^\prime)$ is given by the group structure
of $\hat{\Gamma}$.

It suffices
to define an equivalence of symmetric monoidal functors
$\Psi\colon E^\theta\rightarrow E^\epsilon$.  This is defined by
$\Psi_\gamma|a|:= \gamma^*\hat q(a)$.
This is well-defined as
\begin{align*}
  \Psi_\gamma|a|
  &= \gamma^*\hat q(a)
   = \gamma^*\hat q(b)\gamma^*\hat q(b)^{-1}\gamma^*\hat q(a)
   = \gamma^*\hat q(a-b)\gamma^*\hat q(b)\\
   &= \gamma^*\hat q(\partial e) \gamma^*\hat q(b)
  = \Psi_\gamma(\gamma^*\hat q(\partial e)|b|).
\end{align*}

We now verify that $\Psi$ is a natural transformation.  Let
$g:\Sigma\rightarrow X$ be a cobordism from $\gamma_0$ to $\gamma_1$.
Let
$f\in C_2(\Sigma)$ represent $[\Sigma]\in
H_2(\Sigma,\partial \Sigma)$ and be such that $a_0-a_1=\partial f$,
where $|a_i|\in E^\theta(\gamma_i)$.
Then we have
 \begin{align*}
    \Psi_{\gamma_1}\bigl( E^\theta(g) |a_0|\bigr)
    &= \Psi_{\gamma_1}\bigl( g^*\theta_\epsilon(f)  |a_1|\bigr)
     = g^*\hat q(\partial f)\Psi_{\gamma_1}|a_1|
     =  g^*\hat q(a_0-a_1)\gamma_1^*\hat q(a_1)\\
    &= \gamma_0^*\hat q(a_0) \gamma_1^*\hat q(a_1)^{-1}  \gamma_1^*\hat q(a_1)
     = \Psi_{\gamma_0}|a_0| = E^\epsilon(g)(\Psi_{\gamma_0}|a_0|).
\end{align*}
Finally we verify, that $\Psi$ is compatible with the monoidal structures.
 \begin{align*}
    \Psi_{\gamma_0\sqcup\gamma_1}\bigl(\Phi^\theta_{\gamma_0,\gamma_1}
        (|a_0|\otimes|a_1|)\bigr)
    &= \Psi_{\gamma_0\sqcup\gamma_1} (|a_0\sqcup a_1|)
     = (\gamma_0\sqcup \gamma_1)^*\hat q (a_0+a_1)\\
    &= \gamma_0^*\hat q (a_0) \gamma_1^*\hat q (a_1)
     = \Phi^\epsilon_{\gamma_0,\gamma_1} ( \gamma_0^*\hat q (a_0)\otimes
                        \gamma_1^*\hat q (a_1))\\
    &= \Phi^\epsilon_{\gamma_0,\gamma_1}\bigl(\Psi_{\gamma_0}|a_0|
         \otimes \Psi_{\gamma_1}|a_1|\bigr).
\end{align*}

Thus $E^\theta$ and $E^\epsilon$ are isomorphic and the left-hand
square of the diagram commutes.  As mentioned above, the Five Lemma
can now be invoked to prove that $\tau$ is an isomorphism.
\end{proof}


\section{The line bundle on loop space}
\label{Section:LoopSpaceBundle}

Associated to a gerbe with connection on a manifold $X$ is a line
bundle with connection on $\cL X$, the free loop space on $X$ thought
of as an infinite dimensional manifold (see
\cite[Chapter 6]{Brylinski:LoopSpacesBook}).  Recalling from the
introduction that the second Cheeger-Simons group classifies line
bundles with connection, this association can be viewed as the
transgression map $t\colon
\widehat H^3(X)\to \widehat H^2(\cL X)$ which is described below.
Alternatively, given a thin-invariant field theory $E\colon
\cT_X\to\vect_1$ we can restrict this to a functor on the path
category of the loop space, which gives us a line bundle with
connection on the free
loop space.  Not surprisingly, since this is where the
origins of the definition of a thin-invariant field theory lie, these
two ways of getting a line bundle coincide.

\begin{thm}
    \setcounter{test}{\value{thm}}\addtocounter{test}{-1}
    \refstepcounter{test}\label{thmnumber:LoopSpaceBundle}
  A thin-invariant field theory can be restricted to the path category
  of the free loop space giving a line bundle with connection on the
  free loop space. This line bundle with connection is isomorphic to
  the transgression of the gerbe associated to the thin-invariant
  field theory.
\end{thm}
\begin{proof}
  First we need to describe the transgression map $t\colon \widehat
  H^3(X)\to \widehat H^2(\cL X)$.  Given $S\in \widehat H^3(X)$ a gerbe
  holonomy, define $t(S)\in \widehat H^2(\cL X)$ as follows.  If $\gamma\colon
  S_n\to \cL X$ is a smooth map, then we have an induced map $\check
  \gamma\colon S^1\times S_n\rightarrow X$ given by $\check
  \gamma(r,s):=\gamma(s)(r)$: set $t(S)(\gamma):=S(\check \gamma)$.
  (The curvature of $S$ is transgressed as $t(c)=\pi_*\mathrm{ev}^*c$
  where $\mathrm{ev}\colon S^1\times \cL X\rightarrow X$ is the evaluation
  map, and $\pi_*\colon\Omega^3(S^1\times \mathcal
  LX)\to\Omega^2(\mathcal LX)$ is integration over the fibre.)

  Now if $E$ is a thin-invariant field theory with holonomy $S_E$ to
  verify that the restriction $E^\prime\colon P\mathcal L X\to \vect_1$ is a
  line bundle coinciding with the transgression $t(S_E)$ it suffices
  to compare holonomies.  Let $\gamma \colon S^1\rightarrow \cL X$ be
  a smooth loop. We view $\gamma $ as a map $\check{\gamma }\colon
  S^1\times S^1\rightarrow X$, $\check{\gamma }(r,s):=\gamma (s)(r)$.
  On the one hand the holonomy of $E^\prime$ along $\gamma $ is then given by
  $E(\check{\gamma })\in\mathrm{Aut}(\bC)\cong \bC^\times$.  On the
  other hand, the holonomy of $t(S_E)$ along $\gamma $ is (by the
  definition of $S_E$) equal to $S_E(\check{\gamma })=E(\check
  \gamma)$.
\end{proof}

\sh{Acknowledgements}
The second author would like to thank Roger
Picken for comments on an earlier version. The third author was
supported by a grant from the EPSRC.


\appendix

\section{Appendix: The Cheeger-Simons group}
\label{ap:csg}

The definition of the Cheeger-Simons group given in the main body of
the text is non-standard.  Our definition is based on manifolds and
maps to $X$ whereas the original definition of Cheeger and Simons uses
chains in $X$. In this appendix we prove that these two definitions
are equivalent.

First recall the usual definition of the Cheeger-Simons group. Let
$Z_2X$ be the group of smooth two-cycles in $X$. A {\em
differential character} is a pair $(f,c)$ where $f$ is a
homomorphism $f\colon Z_2X \ra \bC^\times$ and $c$ is a closed
three-form such that if $B$ is a three-chain then

\begin{equation}
\label{Eqn:CS3formcondition}
  f(\partial B) = \exp  \left(2\pi i \int_B c\right).
\end{equation}
The collection of differential characters forms a group which we will
denote by $\widehat{\mathbb H}^3(X)$.  This is the usual definition of
the Cheeger-Simons group, though the index ``$3$'' is\/ {\it \'a la\/}
Brylinski as opposed to the ``$2$'' used by Cheeger and Simons.

In this paper we considered the group $\widehat H^3(X)$ consisting of
pairs
$(S,c)$ where $S$ is a
$\bC^\times$-valued function on the space of maps of closed
surfaces to $X$, and $c$ is a three-form such that if $v\colon V\ra X$ is
a map of a three-manifold to $X$ then
\begin{equation}
\label{Eqn:Our3formcondition}
  S(\partial v) = \exp \left(2\pi i \int_V v^*c\right).
\end{equation}
We can now show that these two groups are essentially the same.

\begin{thm}
  For a smooth manifold $X$, the group $\widehat H^3(X)$ is
  canonically isomorphic to the Cheeger-Simons group $\widehat
  {\mathbb{H}}^3(X)$.
\end{thm}

\proof
Define a function $\widehat H^3(X)\ra\widehat {\mathbb{H}}^3(X)$
via $(S,c)\mapsto(f_S,c)$ where $f_S$ is defined as follows.  By
the isomorphism from $\MSO_2(X)$ to $H_2(X,\bZ)$, for a smooth
two-cycle $y$, there is a closed smooth $X$-surface $g\colon
\Sigma \ra X$ with fundamental cycle $d\in H_2(\Sigma,\bZ)$ such
that $[g_*(d)]=[y]$ in $H_2(X,\bZ)$, ie.\ there is a smooth
three-chain $B$ such that $\partial B = -g_*(d) + y$.  Now define
\[
  f_S(y) = S(g)\exp\left(2\pi i \int_B c\right).
\]

First we must show that this is well-defined, ie.\ that it is
independent of the choices made and that $(f_S,c)$ satisfies
(\ref{Eqn:CS3formcondition}).  Suppose we are given $g^\prime$,
$d^\prime$ and $B^\prime$ such that $\partial B^\prime -g^\prime_*(d^\prime) + y$ then
\[
  g_*(d) - g^\prime_*(d^\prime)
    = ( - \partial B+ y) - ( - \partial B^\prime + y)
    = \partial (B^\prime - B)
\]
showing that  $[g_*(d)] =  [g^\prime_*(d^\prime)]\in H^2(X,\bZ)$ and
so $[g]=[g']\in MSO_2(X)$, meaning that there is an $X$-three-manifold
$v\colon V \ra X$ such that $\partial v g-g^\prime$. Now observe that
\[
  S(g)/S(g^\prime)
     = S(g-g^\prime)
     = S(\partial v)
     = \exp \left(2 \pi i \int_V v^*c\right).
\]
Choosing a  fundamental cycle (relative to the boundary) $D$ for $V$ such that
$\partial D =  d - d^\prime$ we get that $v_*(D)+B-B'$ is a cycle, and so
$\exp\left( 2 \pi i \int_{D+B - B^\prime} c\right) = 1$.  Hence
\begin{align*}
  S(g)\exp\left(2 \pi i \int_B c\right)
     &= S(g^\prime)\exp \left(2 \pi i \int_D c\right)
                  \exp\left(2 \pi i \int_B c\right)\\
     &= S(g^\prime)\exp \left(2 \pi i \int_{B'} c\right),
\end{align*}
so that $f_S$ is independent of the choices made.

We must also show that $(f_S,c)$ satisfies
(\ref{Eqn:CS3formcondition}).  If $B$ is a smooth three-chain, then to
apply $f_S$ to $\partial B$, we may choose $g$ and $d$ above to be
trivial, so we see immediately that $f_S(\partial B) = \exp \left(2
\pi i \int_B c\right)$ as required.


We have a well defined map going one way, so we wish to define a map
going the other way, $ \widehat {\mathbb{H}}^3(X)\ra\widehat H^3(X)$,
which we do via $(f,c)\mapsto( S_f,c)$ where $S_f(g\colon\Sigma\to
X):=f(g_*(d))$, for $d$ a fundamental two-cycle for $\Sigma$.  We will
show that this is an inverse to the above map.  We must first show
this is a well-defined. Suppose we are given another fundamental
two-cycle $d^\prime$, then we can find a three-cycle $e$ in $\Sigma$
such that $\partial e = d-d^\prime$. Observe that
\begin{align*}
  f(g_*(d))/f(g_*(d^\prime))
    &= f(g_*(\partial e))
     = f(\partial(g_*(e)))\\
    &= \exp\left(2 \pi i \int_{g_*(e)} c\right)
     =  \exp\left(2 \pi i \int_{\Sigma} g^*c\right).
\end{align*}
However, $\int_{\Sigma} g^*c= 0$ since $c$ is a three-form and
$\Sigma$ a surface. This shows that
$f\left(g_*\left(d\right)\right)= f(g_*(d^\prime))$.

Now we will verify that $(S_f,c)$ satisfies
(\ref{Eqn:Our3formcondition}).  Given $v\colon V\ra X$, choose a
(relative) fundamental cycle $B$ such that $\partial B$ is a
fundamental cycle for $\partial V$, then
\[
   S(\partial v) = f\left(\partial v_* (\partial B)\right)
                 = f(\partial (v_*B))
                 = \exp\left(2 \pi i \int_{v_*B} c\right)
                 = \exp\left(2 \pi i \int_V v^*c\right).
\]
This shows that $(S_f,c)$ is a well defined element of  $\widehat
H^3(X)$.

Finally we show that the two maps are inverses, ie.\ that $f_{S_f}=f$
and $S_{f_S}=S$.  For the first equality, let $y$ be a smooth
two-cycle in $X$ and choose (as before) a map $g\colon \Sigma \ra X$,
a fundamental two-cycle $d$ and a three-chain $B$ such that $\partial
B = -g_*(d) + y$.
Then
\begin{align*}
  f_{S_f}(y)
    &= S_f(g) \exp\left( 2 \pi i \int_B c\right)
     = f(g_*(d))\exp\left( 2 \pi i\int_B c\right)
     = f(g_*(d))f(\partial B)\\
    &= f(g_*(d)) f(-g_*(d) + y)
     = f(g_*(d)) f(g_*(d))^{-1} f(y)
     = f(y).
\end{align*}
For the second equality, the $B$ can be chosen trivially so that
$$
   S_{f_S}(g) = f_S(g_*(d)) = S(g)\exp\left(2 \pi i \int_B c\right) = S(g).
\eqno{\qed}
$$

The equivalence between the bordism and chain definition of the
Cheeger-Simons group presented in this section is a phenomenon of the
particular low dimension we are working in. For higher dimensions
there is a difference between bordism and homology. It is possible,
however, to define a variant of thin invariant field theory based on chains in $X$
for any dimension $n$ and such theories are related to the
Cheeger-Simons groups in a similar fashion to that
presented in this paper (see \cite{turner}) .


\section{Appendix: Symmetric monoidal categories}\label{ap:defns}
In this appendix we reproduce, for convenient ease of access, the
categorical definitions pertinent to this paper. For further details
see for example \cite{bo:hca2}.

\begin{defn}\label{def:moncat}
  A {\em monoidal category} is a category $\cC$ equipped with a
  bifunctor $\otimes\colon \cC\times\cC \ra \cC$ and an object $\bf1$,
  the unit, together with the following structure isomorphisms:\\
(i)\qua for every triple $A,B,C$ of objects, an isomorphism
\[
a_{A,B,C}\colon (A\otimes B) \otimes C \ra A\otimes (B\otimes C)
\]
(ii)\qua  for every object $A$, isomorphisms
\[
l_A\colon \bf1 \otimes A \ra A \qquad \mbox{ and }\qquad r_a\colon
A\otimes \bf1 \ra A .
\]
The above are subject to the following axioms:\\
(1)\qua  The structure isomorphisms are natural (in all variables).\\
(2)\qua  For each quadruple of objects $A,B,C,D$ the following diagram
commutes.
\[
\xymatrix{((A\otimes B)\otimes C)\otimes D \ar[rr]^{a_{A\otimes
   B,C,D}} \ar[d]_{a_{A,B,C}\otimes 1} & & (A\otimes B) \otimes
   (C\otimes D) \ar[dd]^{a_{A,B, C\otimes D}} \\
(A\otimes (B\otimes C))\otimes D \ar[d]_{A_{A,B\otimes C,D}} & & \\
A\otimes ((B\otimes C)\otimes D) \ar[rr]_{1\otimes a_{B,C,D}} & &
A\otimes (B\otimes (C\otimes D))
}
\]
(3)\qua  For each pair of objects $A,B$ the following diagram commutes.
\[
\xymatrix{(A\otimes \bf1 ) \otimes B \ar[r]^{a_{A,\bf1,B}}
  \ar[dr]_{r_A\otimes 1} & A\otimes (\bf1 \otimes B) \ar[d]^{1\otimes
  l_B} \\
& A\otimes B
}
\]
The category is {\em strict} if the structure isomorphisms are
identities.
\end{defn}

\begin{defn}\label{def:symmoncat}
A {\em symmetric} monoidal category is a monoidal category $\cC$
equipped with natural isomorphisms
\[
s_{A,B}\colon A\otimes B \ra B\otimes A
\]
satisfying the following.\\
(1)\qua  For every triple $A,B,C$ of objects the following diagram
commutes.
\[
\xymatrix{
(A\otimes B)\otimes C \ar[rr]^{s_{A,B}\otimes 1} \ar[d]_{a_{A,B,B}} &
&
(B\otimes A)\otimes C \ar[d]^{a_{B,A,C}} \\
A\otimes (B\otimes C) \ar[d]_{s_{A,B\otimes C}} & & B\otimes (A\otimes
C) \ar[d]^{1\otimes s_{A,C}} \\
(B\otimes C)\otimes A \ar[rr]_{a_{B,C,A}} & & B\otimes (C\otimes A)
}
\]
(2)\qua  For every object $A$ the following diagram commutes.
\[
\xymatrix{
A\otimes \bf1 \ar[r]^{s_{A,\bf1}} \ar[dr]_{r_A} & \bf1 \otimes A
\ar[d]^{l_A} \\
 & A
}
\]
(3)\qua  For every pair $A,B$ of objects the following diagram commutes.
\[
\xymatrix{
A\otimes B \ar[r]^{s_{A,B}} \ar@{=}[dr] & B\otimes A
\ar[d]^{s_{B,A}}\\
 & A\otimes B
}
\]
\end{defn}

\begin{defn}\label{def:monfunc}
Let $\cC$ and $\cD$ be monoidal categories. A {\em monoidal functor}
 is a functor $E\colon \cC \ra \cD$ together with the following
 morphisms in $\cD$:\\
(i)\qua  for each pair $A,B$ of objects in $\cC$ a morphism
\[
\Phi_{A,B}\colon E(A)\otimes E(B) \ra E(A\otimes B),
\]
(ii)\qua  a morphism $\epsilon\colon \bf1_{\cD} \ra E(\bf1_{\cC})$.\\
These must satisfy the following axioms.\\
(1)\qua  The $\Phi_{A,B}$ are natural in both $A$ and $B$.\\
(2)\qua  For each triple $A,B,C$ of objects in $\cC$ the following diagram
commutes.
\[
\xymatrix{
(E(A)\otimes E(B))\otimes E(C) \ar[rr]^{a_{E(A),E(B),E(C)}}
\ar[d]_{\Phi_{A,B}\otimes 1} &
&
E(A)\otimes (E(B) \otimes E(C)) \ar[d]^{1\otimes \Phi_{B,C}} \\
E(A\otimes B) \otimes E(C) \ar[d]_{\Phi_{A\otimes B, C}} & &
E(A) \otimes E(B \otimes C) \ar[d]^{\Phi_{A,B\otimes C}} \\
E((A\otimes B)\otimes C) \ar[rr]_{E(a_{A,B,C})} & & E(A\otimes (B\otimes C))
}
\]
(3)\qua  For each object $A$ of $\cC$ the following two diagrams commute.
\[
\xymatrix{
E(\bf1) \otimes E(A) \ar[r]^{\Phi_{\bf1,A}} & E(\bf1 \otimes A)
\ar[d]^{E(l_A)} \\
\bf1 \otimes E(A) \ar[u]^{\epsilon\otimes 1} \ar[r]_{l_{E(A)}} & E(A)
}
\qquad
\xymatrix{
E(A) \otimes E(\bf1) \ar[r]^{\Phi_{A,\bf1}} & E(A \otimes \bf1 )
\ar[d]^{E(r_A)} \\
 E(A)\otimes \bf1 \ar[u]^{1 \otimes \epsilon} \ar[r]_{r_{E(A)}} & E(A)
}
\]

If the categories $\cC$ and $\cD$ are symmetric monoidal then a {\em
  symmetric monoidal functor} is a monoidal functor as above such that
  for every pair $A,B$ of objects in $\cC$ the following diagram
  commutes.
\[
\xymatrix{
E(A) \otimes E(B) \ar[r]^{s_{E(A),E(B)}} \ar[d]_{\Phi_{A,B}} & E(B)
\otimes E(A) \ar[d]^{\Phi_{B,A}} \\
E(A\otimes B) \ar[r]_{E(s_{A,B})} & E(B\otimes A)
}
\]

\end{defn}


\Addresses\recd


\begin{thebibliography}{12}
\bibitem{Barrett:HolonomyAndPathStructures} J. Barrett, Holonomy and
  path structures in General relativity and Yang-Mills theory,
  Int. J. Theor. Phys., {\bf 30} (1991), no.~9, 1171-1215.
\bibitem{bo:hca2} F. Borceux, {\em Handbook of Categorical Algebra 2,
  Categories and Structures}, Encyclopedia of Mathematics and its
  Applications 52, CUP, Cambridge, 1994.
\bibitem{BrightwellTurner:RepresentationOfHomotopySurfaceCategory}
  M. Brightwell and P. Turner, Representations of the homotopy surface
  category of a simply connected space, J. Knot Theory Ramifications
  {\bf 9} (2000), no.~7, 855--864.
\bibitem{Brylinski:LoopSpacesBook} J.-L. Brylinski, {\it Loop spaces,
  characteristic classes and geometric quantization}, Progr. Math.,
  107, Birkh\"auser Boston, Boston, MA, 1993.
\bibitem{CaetanoPicken:AxiomaticHolonomy} A. Caetano and R. Picken, An
  axiomatic definition of holonomy, Int. J. Math., {\bf 5}
  (1994), no.~6, 835-848.
\bibitem{Chaterjee:Thesis}  D. Chatterjee, On the construction of
  abelian gerbs, PhD thesis (Cambridge) (1998).

\bibitem{cheegersimons83}
J.~Cheeger and J.~Simons,
Differential characters and geometric invariants,
In {\em LNM1167}, pages 50--80. Springer Verlag, 1985.

\bibitem{Freed:QuantumGroups} D. Freed, Quantum groups from path integrals,
  in {\it Particles and fields (Banff, AB, 1994)}, 63--107, Springer,
  New York, 1999.
\bibitem{Freed:HigherAlgebraicStructures} D. Freed, Higher algebraic
  structures and quantization, Comm.\ Math.\ Phys. {\bf 159} (1994),
  no. 2, 343--398.
\bibitem{Gajer:GeometryOfDeligneCohomology} P. Gajer, Geometry of Deligne
  cohomology, Invent. Math. {\bf 127} (1997), no.~1, 155--207.
\bibitem{Hitchin:SpecialLagrangianSubmanifolds} N. Hitchin, Lectures on
  special Lagrangian submanifolds, available from: \arxiv{math.DG/9907034}.
\bibitem{MackaayPicken:HolonomyAndParallelTransport} M. Mackaay and
  R. Picken, Holonomy and parallel transport for abelian gerbes,
  Adv. Math. {\bf 170} (2002), 287-339.
\bibitem{Milnor:hCobordismBook} J. Milnor, {\it Lectures on the
  $h$-cobordism theorem}, Princeton Univ. Press, Princeton, N.J., 1965
\bibitem{Murray:BundleGerbes} M. K. Murray, Bundle gerbes, J. London
  Math. Soc. (2) {\bf 54} (1996), no.~2, 403--416
\bibitem{ro:hqft} G. Rodrigues, Homotopy quantum field theories and
    the homotopy cobordism category in dimension 1+1,  J. Knot Theory
    Ramifications {\bf 12} (2003), no.~3, 287--319.
\bibitem{Segal:EllipticCohomology} G. Segal, Elliptic cohomology
  (after Landweber-Stong, Ochanine, Witten, and others), S\'eminaire
  Bourbaki Exp.\ No.\ 695, Ast\'erisque 161-162 (1988), 187--201.
\bibitem{Segal:TopologicalStructures} G. Segal,  Topological
  structures in string theory, Topological methods in the physical
  sciences (London, 2000). R. Soc. Lond. Philos. Trans. Ser. A
  Math. Phys. Eng. Sci.  359 (2001),  {\bf1784}, 1389--1398.
\bibitem{switzer} R. Switzer, {\it Algebraic Topology - Homology and
  Homotopy}, Springer, Classics in Mathematics, Reprint of the 1975 Edition.
\bibitem{Turaev:HomotopyFieldTheoryInDimension2}  V.~Turaev, Homotopy
  field theory in dimension 2 and group-algebras, preprint 
\arxiv{math.QA/9910010}.
\bibitem{turner} P. Turner, A functorial approach to differential
  characters, Algebraic and Geometric Topology, {\bf 4} (2004), 81-93.
\end{thebibliography}
\end{document}